\documentclass{birkjour}
\usepackage{amsmath, amssymb, amsthm}

\newcommand{\N}{{\mathbf N}}

    \newcommand    {\C}{{\mathbf C}}

\def\<{\langle}
\def\>{\rangle}

\newtheorem{theorem}{Theorem}[section]
\newtheorem{proposition}[theorem]{Proposition}
\newtheorem{lemma}[theorem]{Lemma}
\newtheorem{corollary}[theorem]{Corollary}
\theoremstyle{definition}
\newtheorem{definition}[theorem]{Definition}
\theoremstyle{remark}
\newtheorem{remark}[theorem]{Remark}
\numberwithin{equation}{section}

\newtheorem{problem}{\bf Problem}[section]

\renewcommand{\theequation}{\thesection.\arabic{equation}}
    
\begin{document}

\title[Generic monodromy group]{Generic monodromy group of Riemann surfaces for inverses to entire functions of finite order}

\author[L. Zelenko]
{ Leonid Zelenko} 

\address{%
	Department of Mathematics \\
	University of Haifa  \\
	Haifa 31905  \\
	Israel}
\email{zelenko@math.haifa.ac.il}

\begin{abstract}
We consider the vector space $E_{\rho,p}$ of entire functions of finite order $\rho\in\N$, whose types are not more than $p>0$, endowed with Frechet topology, which is generated by a sequence of weighted norms. We call a function $f\in E_{\rho,p}$ {\it typical} if it is surjective and has an infinite number critical points such that each of them is non-degenerate and all the values of $f$ at these points are pairwise different. We prove that the set of all typical functions contains a set which is $G_\delta$ and dense in $E_{\rho,p}$. Furthermore, we show that inverse to any typical function has Riemann surface whose monodromy group coincides with finitary symmetric group $\mathrm{FSym}(\N)$ of permutations of naturals, which is unsolvable in the following strong sense: it does not have a normal tower of subgroups, whose factor groups are or abelian or finite.  As a consequence from these facts and Topological Galois Theory, we obtain that generically (in the above sense)  for $f\in E_{\rho,p}$ the solution of equation $f(w)=z$ cannot be represented via $z$ and complex constants   by a finite number of the following actions: algebraic operations (i.e., rational ones and solutions of polynomial equations) and quadratures (in particular, superpositions with elementary functions). 
\end{abstract}



\maketitle

\tableofcontents

\section{Introduction}\label{sect:introduct}
\setcounter{equation}{0}

Let $w(z)$ be a multivalued function, defined implicitly as solution of the equation $f(w)=z$, where $f(w)$ be an entire function. It is intuitively clear that in extremely rare cases $w(z)$ can be represented in ``finite terms'', i.e. by a finite number of rational operations and superpositions with elementary 
functions. In this work we try to give a precise meaning to the words ``rare cases''. Notice that in the case where $f(w)$ is a monic polynomial $P_m(w)$, $w(z)$ can be not representable in radicals, as it was shown by Ruffini and Abel \cite{Z}. In Theorem \ref{mainthpol} of the present paper we give a sufficient condition for symmetry of monodromy group of Riemann surface for the inverse polynomial $P_m^{-1}(z)$ in terms of non-vanishing of so called {\it Critical Values  Discriminant} (Definition \ref{defdisccritval}). In particular, for $m\ge 5$ we obtain a sufficient condition in terms of coefficients for 
non-representability in radicals of solutions of the equation $P_m(w)=z$ (Corollary \ref{cornonsolvrad}). This means that $P_m^{-1}(z)$ is not representable in radicals for any point of the space $\C^m$ of coefficients of $P_m(w)$, maybe except of points belonging to a principal algebraic variety, which is the set of zeros of the Critical Values Discriminant. In particular, this means that the set of points in $\C^m$, corresponding to polynomials $P_m$, for which $P_m^{-1}$ is representable in radicals, has zero $2m$-dimensional Lebesgue measure. Let us notice that apparently these  facts can be inferred from results of the classical work of J.F. Ritt \cite{Ritt}, where all the polynomials $P_m$ have been described, for which $P_m^{-1}$ is represented by radicals. 

But in general case  we cannot use the above approach since as far as we know, the concepts of resultant and discriminant are not defined for arbitrary entire functions (except some specific cases - \cite{Goh-Ler}, \cite{Kyt-Nap}). Furthermore, for these functions there is no analogue of the Fundamental Theorem of Algebra and, in particular, they can be non-surjective.
Our approach is that instead of considering a space of sequences of coefficients of  expansions of entire functions in power series, we restrict ourselves on the class $E_{\rho,p}$ of entire functions of finite order $\rho\in\N$, whose types are not more than $p>0$, endowed with Frechet topology, which is generated by a sequence of weighted norms. For a neighborhood $U$ of each point $w_0\in E_{\rho,p}$ and for a disk $D_R\subset\C$ we define the notion of {\it Critical Values Truncated Discriminant} $\mathrm{CVD}_{U,\,R}(w)$ (Definition \ref{dftruncCVD}). To this end we consider the functionals $G(w,z)=w(z)$ and $F(w,z)=w^\prime(z)$ $(w\in E_{\rho,p},\,z\in\C)$ as defined in $E_{\rho,p}\times\C$. This point of view enables us to use methods, applied in theory of functions of several complex variables \cite{Shab}. In particular, Definition \ref{dftruncCVD} of Critical Values Truncated  Discriminant $\mathrm{CVD}_{U,\,R}(w)$ is based on a semi-global version of Weierstrass Preparation Theorem (Proposition \ref{prreprWeir}), which permits us to reduce this definition in a local manner to  a polynomial case. We show that for $w\in U$ $\mathrm{CVD}_{U,\,R}(w)\neq 0$ if and only if the function $w(z)$ has in the disk $D_R$ pairwise different critical values (Proposition \ref{profdiffercritpoints}).  Furthermore, we prove that $\mathrm{CVD}_{U,\,R}(w)$ is holomorphic in $U$ (Proposition \ref{prholofdisc}). To this end we use the theory of analytic functions, defined in a topological vector space (\cite{BS}, \cite{Tay}, \cite{Tay1}). 

We call a function $f\in E_{\rho,p}$ {\it typical} if it is surjective and has an infinite number of critical points such that each of them is non-degenerate and all the values of $f$ at these points are pairwise different. Theorem \ref{thgenfrakG}, a central claim of this paper, states that the set of all typical functions contains a set which is $G_\delta$ and dense in $E_{\rho,p}$. The key moment in its proof is following: the subset of functions $w\in U$, for which  $\mathrm{CVD}_{U,\,R}(w)\neq 0$, is open and dense in $U$. This fact is based on holomorphy of $\mathrm{CVD}_{U,\,R}(w)$ and on uniqueness theorem for analytic functions, defined in a Frechet space (Proposition \ref{pruniq}). Furthermore, we show that inverse to any typical function has Riemann surface whose monodromy group coincides with finitary symmetric group $\mathrm{FSym}(\N)$ of permutations of naturals (Theorem \ref{prmongrfottyp}), which is unsolvable in the following strong sense: it does not have a normal tower of subgroups, whose factor groups are or abelian or finite (Theorem \ref{prcon1}, Definition \ref{defofabfinsolv}).  As a consequence from these facts and Topological Galois Theory, developed by Askold Khovanskii and collaborators (\cite{Kh2}-\cite{Kh1}),  we obtain that generically (in the above sense)  for $f\in E_{\rho,p}$ the solution of equation $f(w)=z$ cannot be represented via $z$ and complex constants   by a finite number of the following actions: algebraic operations (i.e., rational ones and solutions of polynomial equations) and quadratures (in particular, superpositions with elementary functions).
Notice that in the case of entire functions of exponential type ($\rho=1$) it is possible to restrict ourselves on the weighted Banach space $H(\exp(-p|z|))$, defined in Section \ref{subsec:spEp}, since the derivative operator $d/dz$, taking part in the definition of  Critical Values Truncated Discriminant, acts in this space. But for $\rho>1$ it does not act in the space $H(\exp(-p|z|^\rho))$, therefore we are forced to consider the space $E_{\rho,p}$ with a Frechet topology generated by an infinite sequence of weighted norms (see proof of Lemma \ref{lmholofwprim}).

The paper is organized as follows After this Introduction, in Section \ref{sect:prel} we define the notions of Critical Values Discriminant, of Frechet space of entire functions, whose topology is generated by a sequence of weighted norms and of analytic function defined in an open set of topological vector space. In Section \ref{sect:meinres} we formulate main results of the paper, described above. In Section \ref{sect:example} we give an example of typical function. In Section \ref{sect:WeiPrepTh} we formulate and prove a semi-global version of Weierstrass Preparation Theorem, mentioned above. In Section \ref{sect:trcritvaldiscr} we define the Critical Values Truncated Discriminant, state its properties mentioned above and apply it to the proof od local genericity of functions from $E_{\rho,p}$ having pairwise different critical values in a disk (Propositions \ref{prlocgener} and \ref{prdifcritval}). In Section \ref{sect:gensurjfunct} we prove genericity of surjective functions from $E_{\rho,p}$. By the way, Proposition \ref{lmsurjezminepz}, where surjectivity of the function $e^z-\epsilon z\;(\epsilon>0)$ is proved,  plays here an important role. In Section \ref{sect:geninfcritval} we prove genericity of functions from $E_{\rho,p}$, having infinite number of isolated critical points. In Section \ref{sect:proofmainres} we prove the main claims of the paper. In Section \ref{sect: pronlems} we formulate some problems, connected with Theorems \ref{mainthpol} and \ref{thgenfrakG}.
Section \ref{sect:appendix} is Appendix. We prove there some auxiliary claims: in Subsection \ref{subsect:splitinterp} we prove Proposition \ref{lmsplitrots} on splitting of zeros of functions from $E_{\rho,p}$ and Proposition \ref{lmintprob} on specific interpolation polynomials, which permit to construct in the proof of Proposition \ref{prdifcritval} small perturbations of these functions having pairwise different critical values in a disk; in Subsection \ref{subsect:propsomemap} we prove properties of some mappings, used in Sections \ref{sect:gensurjfunct} and \ref{sect:geninfcritval}; in Subsection \ref{subsect:holfuncFrech} we prove some claims on holomorphic functions defined in an open subset of Frechet space and formulate some known claims on them, used in the paper.

\section{Preliminaries}\label{sect:prel}
\setcounter{equation}{0}

\subsection{The case of a polynomial. Critical Values Discriminant}\label{subsect:polyn}

\begin{definition}\label{defcompmatr}
Recall that {\it companion matrix} of a monic polynomial $Q_m(y)=y^m+c_{m-1}y^{m-1}+\dots+c_1y+c_0$ is 
$m\times m$-matrix of the form
\begin{eqnarray*}
C=\left[\begin{array}{lllll}
0&0&\dots&0&-c_0\\
1&0&\dots&0&-c_1\\
0&1&\dots&0&-c_2\\
\cdot&\cdot&\dots&\cdot&\cdot\\
\cdot&\cdot&\dots&\cdot&\cdot\\
\cdot&\cdot&\dots&\cdot&\cdot\\
0&0&\dots&1&-c_{m-1}
\end{array}
\right]
\end{eqnarray*}
and having the property $Q_m(y)=\det(C-y I)$.	
\end{definition}

\begin{definition}\label{defdisccritval}
Consider a monic polynomial $P_m(y)=y^m+a_{m-1}y^{m-1}+\dots+a_1y+a_0$,
where $a_j\in\C\;(j=0,1,\dots m-1)$.  Consider the monic normalization of derivative of this polynomial
$\tilde P_m^\prime(y)=y^{m-1}+\frac{(m-1)}{m}a_{m-1}y^{m-2}+\dots+a_1/m$. Denote by
$Disc(\vec c(\vec a))\;\big(\vec c(\vec a)=\big(a_1/m, a_22/m,\dots, a_{m-1}(m-1)/m\big)\big)$ the discriminant of $\tilde P_m^\prime$ (\cite{Lang}, Chapt V, \S 10) and by
$V^\prime$ denote the discriminant variety of $\tilde P_m^\prime$:
$V^\prime=\{\vec a\in \C^m:\;Disc(\vec c(\vec a))=0\}$. Consider the
$(m-1)\times(m-1)$-matrix $C(\vec a)$, which is the companion
	matrix for the polynomial $\tilde P_m^\prime$, that is
$\tilde P_m^\prime(y)=\mathrm {det}(C(\vec a)-yI)$. Consider the matrix
$P_m(C(\vec a))=C^m(\vec a)+a_{m-1}C^{m-1}(\vec a)+\dots+a_1C(\vec
a)+a_0I$ and its characteristic polynomial $U_{m-1}^{\vec
	a}(y)=\mathrm {det}(P_m(C(\vec a))-yI)$. Denote by $\mathrm{CVD}(\vec a)$
the discriminant of $U_{m-1}^{\vec a}$. We call it the {\it $\mathrm{CV}$  (Critical Values) Discriminant} of
the polynomial $P_m$. 
\end{definition}

\begin{remark}\label{remcritval}
The zeros $y_1,y_2,\dots,y_{m-1}$ of
$\tilde P_m^\prime$ are all the critical points of $P_m$ and by the above construction, all the critical values 
\begin{equation*}
P_m(y_1),P_m(y_2),\dots,P_m(y_{n-1})
\end{equation*}
of $P_m$ are all the zeros of the
polynomial $U_{m-1}^{\vec a}(y)$. 
\end{remark}
Let us notice that $\mathrm{CVD}(\vec a)$ is a polynomial of
$a_0,a_1,\dots,a_{m-1}$ of very large degree (by a coarse estimation
the degree is $m(m-1)^2$).
Also consider the set $V_{cvd}=\{\vec a\in
\C^m:\mathrm{CVD}(\vec a)=0\}$, which we call the {\it $\mathrm{CV}$-discriminant
	variety} of $P_m$ .  It is clear that
$V^\prime\subseteq V_{cvd}$.

\subsection{Spaces $H(v)$ and $E_{\rho,p}$}\label{subsec:spEp}

Consider the wighted Banach space $H(v)$ of entire functions, where $v=v(|z|)$ is a radial weight.  It is defined as the space of entire functions $f(z)$ satisfying the condition
\begin{equation*}
\sup_{z\in\C}v(|z|)|f(z)|<\infty
\end{equation*}
with the norm 
\begin{equation}\label{dfnormv}
\Vert f\Vert_v=\sup_{z\in\C}v(|z|)|f(z)|
\end{equation}
(\cite{Gal}). We assume that the function $v(r)\,(r\ge 0)$ is decreasing and satisfies the condition
\begin{equation*}
\forall\,m\in\N:\quad \sup_{r\ge 0} r^mv(r)<\infty,
\end{equation*}
i.e. $H(v)$ contains all the polynomials.
 
The main object in this paper is the complex vector space $E_{\rho,p}$ of entire functions $f(z)$ of finite order $\rho\in\N$, whose types are not more than $p>0$, i.e. of those ones, for which the condition 
\begin{equation}\label{dfexptyp}
\lim_{R\rightarrow\infty}\sup_{z:\,|z|\ge R}\frac{\ln|f(z)|}{|z|^\rho}\le p
\end{equation}
is satisfied. This space is endowed with Frechet topology, generated by the sequence of norms 
\begin{equation}\label{dfnormn}
\Vert f\Vert_n:=\Vert f\Vert_{v_n},\quad\mathrm{where}\quad v_n(|z|)=\exp\big[-(p+1/n)|z|^\rho\big]\;(n\in\N)
\end{equation}
(\cite{Wol}). We see that for any $n\in\N$ and $w\in E_{\rho,p}$ $\Vert w\Vert_n\le \Vert w\Vert_{n+1}$. Recall
that the sets $U_{n,\epsilon}=\{w\in E_{\rho,p}\,\Vert w\Vert_n<\epsilon\}$ form a base of neighborhoods of zero in this topology,
As it is easy to check, $H(v_{n+1})\subset H(v_n)$ for any $n\in\N$ and
\begin{equation}\label{Episintersec}
E_{\rho,p}=\bigcap_{n\in\N}H(v_n).
\end{equation}
Furthermore, the space $E_{\rho,p}$ is metrizable and  since each $H(v_n)$ is a complete normed space, then $E_{\rho,p}$ is complete with respect to a metric generating its topology, i.e. it is a {\it Frechet space} (\cite{Gel-Shil}, Chapt. I, \S 3). In particular, this means that for $E_{\rho,p}$ {\it Baire Category Theorem }is valid, which claims that a countable intersection of open sets, which are dense in the space $E_{\rho,p}$, is dense in it (\cite{D-Sh}, Chapt. I, \S 6, Theorem 9).
 
Denote by $H^\prime(V_n)$ and $E_{\rho,p}^\prime$ the dual spaces to $H(v_n)$ and $E_{\rho,p}$, i.e. they are the vector spaces of all linear continuous functionals on $H(v_n)$ and $E_{\rho,p}$. It is clear that $H^\prime(v_{n})\subset H^\prime(v_{n+1})$ for any $n\in\N$, if one identifies each $g\in H^\prime(v_{n})$ with $g\vert_{H(v_{n+1})}\in H^\prime(v_{n+1})$,  and 
\begin{equation}\label{Epprimunion}
E_{\rho,p}^\prime=\bigcup_{n\in\N}H^\prime(v_n)
\end{equation}
(\cite{Gel-Shil}, Chapt. I, \S 4). 

\subsection{Definition of analytic functions on complex topological vector space}

In \cite{BS}, Definitions 5.5 and 5.6, the concept of analytic function on an open subset of a complex topological vector space is defined via generalized power series. But it will be more convenient for us to use the following equivalent definition (\cite{BS}, Theorem 6.2):

\begin{definition}\label{dfanalfunct}
	Let $U$ be an open subset of a complex topological vector space $E$. We call a function $f:\,U\rightarrow\C$ {\it analytic} (or {\it holomorphic}), if  it is continuous in $U$ and analytic on  affine lines of the form
	\begin{equation}\label{affline}
	L_{w,\,h}=\bigcup_{\lambda\in\C}\{w_0+\lambda h\}\quad (w_0,\,h\in E).
	\end{equation}
\end{definition}

\section{Main results}\label{sect:meinres}
\setcounter{equation}{0}
Let us return to a polynomial $P_m(y)=y^m+a_{m-1}y^{m-1}+\dots+a_1y+a_0$.
Let $\mathrm{Mon}(\vec
a)\;(\vec a=(a_0,a_1,\dots,a_{m-1}))$ be the monodromy group of the
Riemann surface to the multivalued  function
$P_m^{-1}(x)\;(x\in\C)$.
 As usually, we denote by $\mathrm{S}(m)$ the symmetric
group of $m$ elements. In the case of polynomials the  main result is:
\begin{theorem}\label{mainthpol}
Denote  $\vec a=(a_0,a_1,\dots,a_{m-1})$. If
\begin{equation}\label{maincond}
\mathrm{CVD}(\vec a)\neq 0,
\end{equation}
then 
\vskip2mm

\noindent (i) all the critical points of $P_m$ are non-degenerate and its critical values are pairwise different;
\vskip2mm

\noindent (ii) $\mathrm{Mon}(\vec a)=\mathrm{S}(m)$.
\end{theorem}

Theorem~\ref{mainthpol} and the fact that for $m\ge 5$ the group $\mathrm{S}(m)$ is not solvable (\cite{Kh}, Proposition 5.1)  imply:

\begin{corollary}\label{cornonsolvrad}
If $m\ge 5$ and condition (\ref{maincond}) is satisfied, then the
multi-valued function $P_m^{-1}(x)$ is not expressed by radicals.
\end{corollary}

Let us pass to the case of entire functions from the space $E_{\rho,p}$.  Following to \cite{Z}, we shall define the following concept:

\begin{definition}\label{deftyp}
We call an entire function $f(z)$ {\it typical}, if it is surjective and has infinite number of critical points, such that each of them is non-degenerate and the values of $f$ at these points are pairwise different. We shall denote by $\mathfrak{T}$ the set of all typical functions from $E_{\rho,p}$.
\end{definition}

Recall that a subset of a topological space is called $G_\delta$ set, if it ia a countable intersection of open sets. The following theorem is a central claim of this paper:

\begin{theorem}\label{thgenfrakG}
Functions from $\mathfrak{T}$ are generic in $E_{\rho,p}$ in the sense that there is a subset of $\mathfrak{T}$, which is  $G_\delta$ and dense in $E_{\rho,p}$.
\end{theorem}

Recall that a bijective mapping 
$\Delta:\;\N\rightarrow\N$ is called  {\it finitary permutation}, if it is the identity
outside of a finite set. 
Consider the group $\mathrm{FSym}(\N)$ of {\it all} finitary permutations
of the set $\N$. It is called in the literature {\it finitary symmetric permutation group}.

The following claim is valid:

\begin{theorem}\label{prmongrfottyp}
	For any function $f\in\mathfrak{T}$ the monodromy group of Riemann surface of its inverse $f^{-1}$ coincides with $\mathrm{FSym}(\N)$.
\end{theorem}

Following to \cite{Kh}, \cite{Kh1}, we shall mean the solvability of a group in following weakened sense: 

\begin{definition}\label{defofabfinsolv}
	We call an abstract  group $G$ {\it interspersed-finite solvable}, if there exists 
	a decreasing finite series of subgroups of $G$
	\begin{equation}\label{seqofgr}
	\{e\}= G_r\subset G_{r-1}\subset G_{r-2}\subset\dots\subset G_{1}\subset G_0=G,
	\end{equation}
	such that for each $k\in\{0,\,1,\,2,\,\dots,\,r-1\}$ the subgroup $G_{k+1}$ of $G_k$ is normal and the
	quotient  group $G_k/G_{k+1}$ is either abelian or finite. 
\end{definition}

The following claim is valid:

\begin{theorem}\label{prcon1}
The group $\mathrm{FSym}(\N)$ is interspersed-finite unsolvable.
\end{theorem}

 From Topological Galois Theory (\cite{Kh2}-\cite{Kh1}), we obtain the following claims:
 
 \begin{corollary}\label{corrol}
   For $f\in\mathfrak{T}$ the solution of equation $f(w)=z$ cannot be represented via $z$ and complex constants  by a finite number of the following actions: algebraic operations (i.e., rational ones and solutions of polynomial  equations) and quadratures (in particular, superpositions with elementary functions).	
 \end{corollary}  

 As a consequence from above claims we obtain:

\begin{corollary}\label{corro2}
Generically (in the sense, indicated in Theorem '\ref{thgenfrakG})  for $f\in E_{\rho,p}$ the solution of equation $f(w)=z$ cannot be represented via $z$ and complex constants in the manner, indicated in Corollary \ref{corrol}. 
 \end{corollary}
 
 \section{Example of typical function}\label{sect:example}
\setcounter{equation}{0}

Consider the entire function $f(w)=e^w-w$. By Proposition \ref{lmsurjezminepz}, it is surjective.
Its critical points and critical values have the form
\begin{equation*}
w_k=2\pi ki,\quad z_k=1-2\pi ki.
\end{equation*}
It is easy to check that all the critical points $w_k$ are non-degenerate and the critical values $z_k$ are pairwise different. Hence this function is typical. Therefore by Corollary \ref{corrol} the solution of the equation $e^w-w=z$ cannot be represented via $z$ and complex constants in the manner indicated there.

\section{A semi-global version of Weierstrass Preparation Theorem}\label{sect:WeiPrepTh}
\setcounter{equation}{0}

Denote by $D_R\,(R>0)$ the disk 
\begin{equation*}
D_R:=\{z\in\C:\;|z|<R\}
\end{equation*} 
and by $\bar D_R$ denote its closure.

The following infinite-dimensional and semi-global version of Weierstrass Preparation Theorem  is valid:

\begin{proposition}\label{prreprWeir}
Let $\Omega$ be an open subset of a complex Frechet space $E$ and  $F(w,\,z)$ be a holomorphic function on $\Omega\times \C$. Suppose that for some fixed $w_0\in\Omega$ the entire function $F(w_0,\,z)$ has in a disk $D_R$ $m\ge 0$ zeros $z_1^0,\,z_2^0,\,\dots,\,z_m^0$ (counting their multiplicities) and it does not vanish on $\partial D_R$. Then there is a neighborhood $U\subset\Omega$ of $w_0$ such that the function $F(w,\,z)$ has the representation in $U\times D_R$:
\begin{equation}\label{reprWeir}
F(w,\,z)=P_m(z,\,w)\Phi(w,\,z), 
\end{equation}
where 
\begin{eqnarray}\label{dfpolyn}
&&\hskip-5mmP_m(w,\,z)=z^m+a_{m-1}(w)z^{m-1}+\,\dots,\,+a_0(w),\quad\mathrm{if}\quad m\ge 1,\\
&&\hskip-5mmP_m(w,\,z)\equiv 1,\quad\mathrm{if}\quad m=0,\nonumber
\end{eqnarray}
the functions $a_0(w),\,a_1(w),\,\dots,\,a_{m-1}(w)$ are holomorphic in $U$ and the function $\Phi(w,\,z)$ is holomorphic in the domain $U\times D_R$ and it does not vanish in $U\times \bar D_R$.
\end{proposition}
\begin{proof}
 Our arguments will be similar to ones in Theorem 1 from \cite{Shab}, Chapt. II, \S 8.
In view of continuity of $F(w,\,z)$, compactness of $\partial D_R$ in $\C$ and Rouche's Theorem, there is a neighborhood $U\subset B$ of $w_0$ such that for any fixed $w\in U$ $F(w,\,z)\neq 0$ on $\partial D_R$ and the $F(w,\,z)$ has in $D_R$ $m$ zeros 
\begin{equation*}
z_1(w),\,z_2(w),\,\dots,\,z_m(w)
\end{equation*}
 (counting their multiplicities). Let us define the polynomial \eqref{dfpolyn} as having the same zeros:
\begin{equation*}
P_m(w,z)=\prod_{\nu=1}^m(z-z_\nu(w)).
\end{equation*}
Let us show that the coefficients $a_k(w)$ are holomorphic in $U$. By Proposition \ref{prholint}, for any entire function $\omega(z)$ the representation
\begin{equation*}
\sum_{\nu=1}^n\omega(z_\nu(w))=\frac{1}{2\pi i}\oint_{\partial D_R}\omega(z)\frac{\partial F(w,z)/\partial z}{F(w,z)}\,\mathrm{d}z.
\end{equation*}
is valid, where the integral is holomorphic in $U$.
Taking in particular $\omega(z)=z^m\,(m\in\N)$, we obtain that Newtonian sums $\sum_{\nu=1}^n z_\nu^m(w)$ are holomorphic in $U$. But the well known Newton identities yield representations of $a_k(w)$ as polynomials of Newtonian sums, hence, in view of claim (iii) of Proposition \ref{properholfunc}, they are holomorphic in $U$. Now let us show that 
the function 
\begin{equation*}
\Phi(w,\,z)=\frac{F(w,\,z)}{P_m(w,z)}
\end{equation*}
has the desired properties. We see that for any fixed $w\in U$ it is holomorphic with respect to $z$ in the disk $D_R$ and does not vanish in $\bar D_R$, since it has the zeros of the same order as $P_n(w,z)$ only at the points $z_\nu(w)$. Then for any $w\in U$ it is represented  by Cauchy integral by the variable $z$:
\begin{equation*}
\Phi(w,\,z)=\frac{1}{2\pi i}\oint_{\partial D_R}\frac{F(w,\,\zeta)}{P_m(w,\zeta)}\frac{\mathrm{d}\zeta}{(\zeta-z)}.
\end{equation*}
We see that $\Phi(w,\,z)$ is defined by this representation also at the points $z=z_\nu(w)$. Furthermore, since $P_n(w,z)\neq 0$ on $\partial D_R$ for any $w\in U$, we obtain using all claims of Proposition \ref{properholfunc}, that the integrand in this representation is  holomorphic with respect to $w$ in $U$ for any fixed $\zeta\in \partial D_R$ and continuous in $U\times\partial D_R$. Hence by Proposition \ref{pronintpar}, for each fixed $z\in D_R$ the function $\Phi(w,\,z)$  is holomorphic in $U$. Thus, by  Propositions \ref{prHartogs} (which is a generalization of Hartogs Theorem), $\Phi(w,\,z)$ is holomorphic in $U\times D_R$.
\end{proof} 

\section{Critical Values Truncated Discriminant and local genericity of functions with pairwise different critical values}\label{sect:trcritvaldiscr}
\setcounter{equation}{0}

Before to give definition of the notion, mentioned in the title of this subsection, we need the following claim:
\begin{lemma}\label{lmholofwprim}
\noindent (i) $w^\prime\in E_{\rho,p}$ for any $w\in E_{\rho,p}$;
\vskip2mm

\noindent(ii) The functional 
$F(w,\,z):=w^\prime(z)$, 
considered as defined in the complex Frechet space $E_{\rho,p}\times\C$, is holomorphic there;
\vskip2mm

\noindent(iii) The functional 
$G(w,\,z):=w(z)$, 
considered as defined in the complex Frechet space $E_{\rho,p}\times\C$, is holomorphic there.
\end{lemma}
\begin{proof}
(i)	Suppose that $w\in E_{\rho,p}$. Consider the circle $C_1(z)\subset\C$, whose center is $z$ and radius $1$, Using Cauchy inequality and \eqref{dfnormv}, \eqref{dfnormn}, \eqref{Episintersec}, we obtain for $n\in\N$:
\begin{eqnarray}\label{estwprim}
&&\hskip-10mm|w^\prime(z)|\le\max_{\zeta\in C_1(z)}|w(\zeta)|\le \exp\big((p+1/n)|z|^\rho\big)\times\nonumber\\
&&\hskip-10mm\exp\big((p+1/(n+1))(|z|+1)^\rho-(p+1/n)|z|^\rho\big)\times\\
&&\hskip-10mm\max_{\zeta\in C_1(z)}\exp\big(-(p+1/(n+1))\,|\zeta|^\rho\big)|w(\zeta)|\le
A\exp\big((p+1/n)|z|^\rho\big)\Vert w\Vert_{n+1}, \nonumber 
\end{eqnarray}
where 
\begin{equation*}
A=\sup_{t\ge 0}\Big(\exp\Big(\big(p+\frac{1}{n+1}\big)(t+1)^\rho-\big(p+\frac{1}{n}\big)t^\rho\Big)\Big)<\infty,
\end{equation*}
since the function under exponent sign is continuous and it is negative for sufficiently large $t$. This means that  $w^\prime\in H(v_n)$ or any $n\in\N$. Then in view of \eqref{Episintersec},  $w^\prime\in E_{\rho,p}$. Claim (i) is proven.
	
(ii) By Proposition \ref{prHartogs}, we should prove that $F(w,\,z)$ is holomorphic by $z$ for any fixed $w\in E_{\rho,p}$ and by $w$ for any fixed $z\in\C$. The first claim is obvious. Let us prove the second one. In view of Definition \ref{dfanalfunct}, we should show that for any fixed $z\in\C$ the function $F(w,\,z)$ is continuous in $E_{\rho,p}$ and analytic on any affine line of the form \eqref{affline}.
 We obtain that 
\begin{equation*}
F(w,z)\vert_{w=w_0+\lambda h}=w_0^\prime(z)+\lambda h^\prime(z),
\end{equation*}
i.e. the last function is holomorphic with respect to $\lambda$.
Let us show  that for any fixed $z\in C$ $F(w,\,z)$ is continuous, i.e. the linear functional $w^\prime(z)$ belongs to 
the dual space $E_{\rho,p}^\prime$. We see from estimate \eqref{estwprim} that for fixed $z\in\C$ and $n\in\N$ the functional $w^\prime(z)$ belongs to $H^\prime(v_{n+1})$, hence in view of \eqref{Epprimunion}, it belongs to $E_{\rho,p}^\prime$. Claim (ii) is proven.

The proof of claim (iii) is similar to the proof of claim (ii) and even simpler than the latter.
\end{proof}

We turn now to the promised definition. 

\begin{definition}\label{dftruncCVD}
Suppose that for some $w_0\in E_{\rho,p}$ the function $w_0^\prime (z)$ has $m\ge 0$ zeros in a disk $D_R$ and $w_0^\prime (z)\neq 0$ on $\partial D_R$. Then by Proposition \ref{prreprWeir} and Lemma \ref{lmholofwprim}, there is a neighborhood $U$ of $w_0$ in $E_{\rho,p}$ such that  the function $F(w,\,z):=w^\prime(z)$ has the representation \eqref{reprWeir}-\eqref{dfpolyn} in $U\times D_R$, where the coefficients $a_0(w),\,a_1(w),\,\dots,\,a_{m-1}(w)$ of the polynomial $P_m(w,\,z)$ are holomorphic in $U$. In the case where $m\ge 1$ let $C(w)$ be the companion matrix for the polynomial $P_m(w,\,z)$ (Definition \ref{defcompmatr}) and the matrix $W(C)(w)$ be the function $w$ of this matrix defined in the usual manner
\begin{equation}\label{funcofcompmat}
W(C)(w);=-\frac{1}{2\pi i}\oint_{\partial D_R}w(\zeta)\big(C(w)-\zeta I\big)^{-1}\,\mathrm{d}\zeta
\end{equation}
(\cite{D-Sh}, Chapt. VII, \S 3, Definition 9). Consider the characteristic polynomial $D_{m,\,W}(\lambda,\,w)=\det\big(W(C)(w)-\lambda I\big)$ of the matrix $W(C)(w)$. We call its discriminant $\mathrm{Disc}\big(D_{m,\,W}(\lambda,\,w)\big)$  {\it  $\mathrm{CV}$ (Critical Values) $R$-Truncated Discriminant} for a function $w\in U$ on the disk $D_R$ and denote it by $\mathrm{CVD}_{U,\,R}(w)$. In the case where $m=0$ we put $\mathrm{CVD}_{U,\,R}(w)=1$ for any $w\in U$.
\end{definition}

The following properties of $\mathrm{CVD}_{U,\,R}(w)$ are valid:

\begin{proposition}\label{profdiffercritpoints}
(i)  $\mathrm{CVD}_{U,\,R}(w)\neq 0$ for some $w\in U$ if and only if all the critical values of the function $w(z)$ in $D_R$ are pairwise different;
\vskip2mm

\noindent (ii) If $\mathrm{CVD}_{U,\,R}(w)\neq 0$, all the critical points of $w(z)$ in $D_R$ are non-degenerate.
\end{proposition}
\begin{proof}
(i) According to Definition \ref{dftruncCVD}, $\mathrm{CVD}_{U,\,R}(w)\neq 0$ if and only if all the 
eigenvalues of the matrix $W(C)(w)$, lying in $D_R$, have multiplicity one, i.e. they are pairwise different (\cite{Lang}, Chapt. 5, \S 10, Corollary of Proposition 4 and Proposition 5). In view of definition \eqref{funcofcompmat} and Spectral Mapping Theorem (\cite{D-Sh}, Chapt. VII, \S 3, Theorem 11), these eigenvalues coincide with $w(z_k)\,(k=1,2,\dots,m)$, where $z_k$ are all the eigenvalues of the matrix $C(w)$, lying in $D_R$, which in turn are zeros of the derivative $w^\prime(z)$ (i.e., they are critical points of $w(z)$ in $D_R$). Claim (i) is proven.

(i) By claim (i), if $\mathrm{CVD}_{U,\,R}(w)\neq 0$, then the critical values $w(z_k)\,(k=1,2,\dots,m)$ of the function $w(z)$ in $D_R$ are pairwise different. hence all the critical points $z_k\,(k=1,2,\dots,m)$ of $w(z)$ in $D_R$ are pairwise different, i.e. each of them is non-degenerate. Claim (ii) is proven.
\end{proof}

\begin{proposition}\label{prholofdisc}
The $CV$ Truncated Discriminant $\mathrm{CVD}_{U,\,R}(w)$ is holomorphic in $U$. 
\end{proposition}
\begin{proof}
Consider $G(w,\,z):=w(z)$ as defined in the complex Frechet space $E_{\rho,p}\times\C$ and write \eqref{funcofcompmat} in the form:
\begin{equation*}
W(C)(w);=-\frac{1}{2\pi i}\oint_{\partial D_R}G(w, \zeta)\big(C(w)-\zeta I\big)^{-1}\,\mathrm{d}\zeta
\end{equation*}
Since by Proposition \ref{prreprWeir} and Lemma \ref{lmholofwprim} the coefficients
\begin{equation*}
a_0(w),\,a_1(w),\,\dots,\,a_{m-1}(w)
\end{equation*}
  of the polynomial $P_m(w,\,z)$ are holomorphic in $U$, then all the enters of its companion matrix $C(w)$ are holomorphic in $U$ (see Definition \ref{defcompmatr}). Furthermore since the circle $\partial D_R$ does not contain the zeros of the polynomial $P_m(w,\,z)$ (which are eigenvalues of the matrix $C(w)$), then $\det(C(w)-\zeta I)\neq 0$ for any $\zeta\in\partial D_R$. On the other hand, for $\zeta\in\partial D_R$ the enters of the matrix $\big(C(w)-\zeta I\big)^{-1}$ have the form:
 \begin{equation*}
 \frac{A_{j,k}(w,\zeta)}{\det(C(w)-\zeta I)}, 
 \end{equation*}
where $A_{j,k}(w,\zeta)$ and $\det(C(w)-\zeta I)$ are polynomials of 
\begin{equation*}
a_0(w),\,a_1(w),\,\dots,\,a_{m-1}(w),\,\zeta.
\end{equation*}
Furthermore, by claim (iii) of Lemma \eqref{lmholofwprim}, the function $G(w,z)$ is holomorphic in $U\times\C$. Hence by Propositions \ref{pronintpar} and \ref{properholfunc}, the enters of the matrix $W(C)(w)$ are holomorphic in $U$. But since the discriminant of the polynomial $\det\big(W(C)(w)-\lambda I\big)$ is a polynomial of the entries of $W(C)(w)$, then using claim (iii) of Proposition \ref{properholfunc}, we obtain that $\mathrm{CVD}_{U,\,R}(w)$ is holomorphic in $U$.
\end{proof}

Following to \cite{Ad}, we introduce the concept:

\begin{definition}\label{dfprincvar}
A subset $\mathcal{V}$ of an open set $U$ in a topological vector space $B$ is called {\it principal analytic variety}, if there is a holomorphic functional $\Phi$ on $U$, such that it does not vanish identically on any connected component of $U$ and $\mathcal{V}=\{w\in U:\;\Phi(w)=0 \}$.
\end{definition}

The main claim of this subsection is following:

\begin{proposition}\label{prlocgener}
Assume that $w_0\in E_{\rho,p}$ and choose $R>0$ so that $w_0(z)$ has no critical points on $\partial D_R$. Then there are a connected neighborhood $U$ of $w_0$ in $E_{\rho,p}$ and a principal analytic variety $\mathcal{V}\subset U$ such that the set of all functions $w\in U$ having only non-degenerate critical points in $D_R$ and taking at them pairwise different values, coincides with $U\setminus\mathcal{V}$.
\end{proposition}
\begin{proof}
Let $U$ be a neighborhood of $w_0$ in $E_{\rho,p}$, mentioned in Definition \ref{dftruncCVD}. Since the space $E_{\rho,p}$ is locally convex, we can choose $U$ such that it is connected. Then by Proposition \ref{profdiffercritpoints}, $\mathrm{CVD}_{U,\,R}(w)\neq 0$ for a function $w\in U$ if and only if all its critical values, corresponding to its critical points in $D_R$, are pairwise different (hence all these critical points are non-degenerate). On the other hand, by Proposition \ref{prdifcritval} there exists a function $\tilde w\in U$ having the last property, i.e. $\mathrm{CVD}_{U,\,R}(\tilde w)\neq 0$. Hence, since by Proposition \ref{prholofdisc} the functional $\mathrm{CVD}_{U,\,R}$ is holomorphic in $U$, the set 
\begin{equation*}
\mathcal{V}=\{w\in U:\;\mathrm{CVD}_{U,\,R}(w)= 0\}
\end{equation*}
is a principal analytic variety in $U$\footnote{In view of Definition \ref{dftruncCVD}, if $w_0(z)$ has no critical points in $\bar D_R$, then $\mathcal{V}=\emptyset$}.
\end{proof}

In the proof of Proposition \ref{prlocgener} we have used the following claim:

\begin{proposition}\label{prdifcritval}
	Let $w_0(z)$ be a function from $E_{\rho,p}$ such that it has in the disk $D_R$ $m$ critical points $z_1^0,\,z_2^0,\,\dots,\,z_m^0$ (counting their multiplicities) and it has no critical points on $\partial D_R$. Then in any open neighborhood $U\subset E_{\rho,p}$ of $w_0$ there is a function $w_{_U}$ such that it has in $D_R$ only $m$ non-degenerate critical points $z_{_{1,U}},\,z_{_{2,U}},\,\dots,\,z_{_{m,U}}$ with pairwise different values 
	\begin{equation*}
	w_{_U}(z_{_{1,U}}),\,w_{_U}(z_{_{2,U}}),\,\dots,\,w_{_U}(z_{_{m,U}})
	\end{equation*}
	 and it has no critical points on $\partial D_R$.
\end{proposition}
\begin{proof}
By the definition of topology on $E_{\rho,p}$, it is sufficient to show tat for any $n\in\N$ and $\epsilon>0$ there is a function $w_{\epsilon,n}\in E_{\rho,p}$ such that $\Vert w_{\epsilon,n}-w_0 \Vert_n<\epsilon$ and it satisfies all the conditions of our claim. Recall that the norms $\Vert\cdot\Vert_n$ and the weights $v_n $are defined by \eqref{dfnormn}. Notice that in view of \eqref{Episintersec}, $w_o\in H(v_n)$ for any $n\in\N$.
By claim (i) of Lemma \ref{lmholofwprim}, $w_0^\prime\in E_{\rho,p}$, hence $w_0^\prime\in H(v_n)$ for any $n\in\N$. By Proposition \ref{lmsplitrots}, applied to the space $H(v_n)$ and the function $w_0^\prime(z)$, for any $\epsilon>0$ there is an entire function  $\theta_{\epsilon,n}(z)$ such that $\theta_{\epsilon}\in H(\tilde v_n)$ with the weight $\tilde v_n(r)=(1+r^2)v_n(r)$,
\begin{equation}\label{estnormvtheta}
\Vert\theta_{\epsilon,n}\Vert_{\tilde v_n}<K\epsilon
\end{equation}
($K>0$ will be chosen below), the function $w_0^\prime(z)+\theta_{\epsilon,n}(z)$ has in $D_R$ $m$ pairwise different zeros $z_1^{\epsilon,n},\,z_2^{\epsilon,n},\,\dots,\,z_m^{\epsilon,n}$ and it does not vanish on $\partial D_R$. Then these zeros are non-degenerate critical points of the function $\tilde w_{\epsilon,n}(z)=w_0(z)+\gamma_{\epsilon,n}(z)$, where $\gamma_{\epsilon,n}(z)=\int_0^{z}\theta_{\epsilon,n}(\zeta)\,\mathrm{d}\zeta$, and $\tilde w_{\epsilon,n}(z)$  has no critical points on $\partial D_R$. Let us show that the function $\gamma_{\epsilon,n}(z)$ belongs to the space $E_{\rho,p}$ and estimate the norms $\Vert\gamma_{\epsilon,n}\Vert_n$. Consider the straight segment, connecting the points $0$ and $z$: $L_z=\{z\in\C: z=t\exp(i\varphi),\, t\in[0,\,|z|] ,\,\varphi=\mathrm{arg}(z)$. We have, using Schwartz inequality:
\begin{eqnarray*}
&&|\gamma_{\epsilon,n}(z)|=\big|\int_{L_z}\theta_{\epsilon,n}(\zeta)\,\mathrm{d}\zeta\big|\le\int_0^{|z|}(1+t^2)
|\theta_{\epsilon,n}(t\exp(i\varphi))|\frac{1}{1+t^2}\,\mathrm{d}t\le\\
&&\Vert\theta_{\epsilon,n}\Vert_{\tilde v_n}\,\int_0^{|z|}\exp((p+1/n)t^\rho)\frac{1}{1+t^2}\,\mathrm{d}t\le\\
&&\Vert\theta_{\epsilon,n}\Vert_{\tilde v_n}\,\sqrt{\int_0^{|z|}\exp(2(p+1/n)t^\rho)\,\mathrm{d}t}\sqrt{\int_0^\infty\frac{\mathrm{d}t}{(1+t^2)^2}}\le\\ 
&&B\Vert\theta_{\epsilon,n}\Vert_{\tilde v_n}\sqrt{\int_0^\infty\frac{\mathrm{d}t}{(1+t^2)^2}}\exp((p+1/n)|z|^\rho),
\end{eqnarray*}
where 
\begin{equation*}
B=\sup_{s\ge 0}\Big(\frac{\sqrt{\int_0^{s}\exp(2(p+1/n)t^\rho)\,\mathrm{d}t}}{\exp((p+1/n)s^\rho)}\Big)<\infty,
\end{equation*}
since the function under brackets is continuous and by L'Hospital's rule,
\begin{equation*}
\lim_{s\rightarrow\infty}\frac{\int_0^{s}\exp(2(p+1/n)t^\rho)\,\mathrm{d}t}{\exp(2(p+1/n)s^\rho)}=\lim_{s\rightarrow\infty}\frac{1}{\rho s^{\rho-1}}<\infty.
\end{equation*}
This means that $\gamma_{\epsilon,n}\in E_{\rho,p}$ and if we shall take in \eqref{estnormvtheta}
\begin{equation*}
K=\frac{1}{2B}\Big(\int_0^\infty\frac{\mathrm{d}t}{(1+t^2)^2}\Big)^{-1/2},
\end{equation*}
then $\Vert\gamma_{\epsilon,n}\Vert_n<\epsilon/2$. Now we shall use Proposition \ref{lmintprob} for choice of a small perturbation of the function $\tilde w_{\epsilon,n}(z)$ in order all critical values will be pairwise different. We choose it in the form $w_{\epsilon,n}(z)=\tilde w_{\epsilon,n}(z)+P_{2m}(z)$, where $P_{2m}(z)$ is the polynomial defined by \eqref{dffoundpol}-\eqref{dfintpol} with $z_k=z_k^{\epsilon,n}\, (k=1,2,\dots,m)$. By Proposition \ref{lmintprob}, $P_{2m}(z)$ satisfies the conditions
\begin{equation*}
P_{2m}(z_j^{\epsilon,n})=y_j,\quad P_{2m}^\prime(z_j^{\epsilon,n})=0\quad (j=1,2,\dots,m). 
\end{equation*}
Hence $w_{\epsilon,n}(z)$ has the same critical points $z_k^{\epsilon,n}$ and $w_{\epsilon,n}(z_k^{\epsilon,n})=\tilde w_{\epsilon,n}(z_k^{\epsilon,n})+y_k\,(k=1,2,\dots, m)$. Since the polynomials $M_{2m}^{j}(z)$ belong to the space $H(v_n)$, we see from 
\eqref{dfintpol} that if $y_j$ are small enough, then $\Vert w_{\epsilon,n}-\tilde w_{\epsilon,n}\Vert_n<\epsilon/2$, hence $\Vert w_{\epsilon,n}-w_0\Vert_n<\epsilon$. On the other hand, it is clear that we can choose these small $y_j$ in order to split the repeating values in the sequence $\{\tilde w_{\epsilon,n}(z_k^{\epsilon,n})\}_{k=1}^n$ in such a manner that all the numbers $w_{\epsilon,n}(z_k^{\epsilon,n})$ will differ from each other.
\end{proof}

\section{Genericity of surjective functions from $E_{\rho,p}$}\label{sect:gensurjfunct}
\setcounter{equation}{0}

Denote by $D_r ^m$ the polydisk  $D_r^m=\times_{k=1}^{m} D_r\;(r>0)$.

The main claim of this subsection is following:

\begin{proposition}\label{prgensurj}
The set $\mathcal{S}$ of all surjective functions from $E_{\rho,p}$ 
is $G_\delta$ and dense in $E_{\rho,p}$.
\end{proposition}
\begin{proof}
By claim (i) of Proposition \ref{lmpropfamexp}, any  non-surjective function from $E_{\rho,p}$ has the form 
\begin{eqnarray}\label{dffalpha}
&&f_{\vec\alpha,\,a,\,b}(z)=b\exp\big(\sum_{k=1}^\rho \alpha_kz^k\big)+a,\\
&&\mathrm{where}\quad |\alpha_\rho|\le p,\;\alpha_k\in\C\;\;\mathrm{for}\;\; k=1,2,\dots,\,\rho-1,\, a, b\in\C\nonumber
\end{eqnarray}
and $\vec\alpha=(\alpha_1,\alpha_2,\dots,\alpha_p)$. 
We see that $f_{\vec\alpha,\,a,\,b}(z)$ is not constant if and only if $\vec\alpha\neq \vec 0$ and $b\neq 0$. Hence
\begin{equation*}
\mathcal{S}^c=E_{\rho,p}\setminus\mathcal{S}=\{f_{\vec\alpha,\,a,\,b}(z)\}_{\vec\alpha\in(\C^{\rho-1}\times\bar D_p) \setminus\{\vec 0\},\,a\in\C,\,b\in\C\setminus\{0\}}\cup\mathfrak{C},
\end{equation*}
where $\mathfrak{C}$ is the set of all constant functions. By Lemma \ref{lmfrakCc}, it is closed in $E_{\rho,p}$.
In order to prove that $\mathcal{S}$ is $G_\delta$ in $E_{\rho,p}$, it is sufficient to show that $\mathcal{S}^c$ is a countable union of closed sets. Denote 
\begin{equation*}
\hskip-10mm\mathcal{K}_{l}=\{f_{\vec\alpha,\,a,\,b}(z)\}_{\vec\alpha\in(\bar D^{\rho-1}_l\times\bar D_p) \setminus D_{\delta_l}^\rho,\,a\in\bar D_l,\,b\in\bar D_l\setminus D_{\delta_l}}\cup\mathfrak{C}.
\end{equation*}
where $l\in\N$ and $\delta_l=p/(l+1)$. Then
\begin{equation}\label{reprcalSpc}
\mathcal{S}^c=\Big(\bigcup_{l\in\N}\mathcal{K}_{l}\Big)\cup\mathfrak{C},
\end{equation}
 In view of claim (ii) of Proposition \ref{lmpropfamexp}, the sets $\mathcal{K}_{l}$ are compact in $E_{\rho,p}$ as continuous images of compact sets $\big((\bar D^{\rho-1}_l\times\bar D_p) \setminus D_{\delta_l}^\rho\big)\times\bar D_l\times\big(\bar D_l\setminus D_{\delta_l}\big)$, hence they are closed in $E_{\rho,p}$. Therefore, in view of \eqref{reprcalSpc}, the set $\mathcal{S}$ is $G_\delta$ in $E_{\rho,p}$.
 
Let us show that $\mathcal{S}$ is dense in $E_{\rho,p}$. To this end it is sufficient to show that for any $n\in\N$, $\sigma>0$ and $f_{\vec\alpha,\,a,\,b}\in\mathcal{S}^c$ there is $g_\sigma\in\mathcal{S}$ such that $\Vert f_{\vec\alpha,\,a,\,b}-g_\sigma\Vert_n<\sigma$. Let us take
\begin{equation}\label{reprgsigma}
 g_\sigma(z)=f_{\vec\alpha,\,a,\,b}(z)-\epsilon\sum_{k=1}^\rho \alpha_kz^k=f_\epsilon\big(\sum_{k=1}^\rho \alpha_kz^k\big) , .
\end{equation}
where $\epsilon>0$ and $f_\epsilon(z)=e^z-\epsilon z$. We see that in the case where
 $\vec\alpha\neq\vec 0$ the polynomial $Q_\rho(z)=\sum_{k=1}^\rho \alpha_kz^k$ is surjective (by  Fundamental Theorem of Algebra). On the other hand, by Proposition \ref{lmsurjezminepz} the function $f_\epsilon(z)$ is surjective. Hence in view of \eqref{reprgsigma}, the function $g_\sigma(z)$ is surjective too.
Thus, since $Q_\rho\in E_{\rho,p}$ for $\vec\alpha\neq 0$ and sufficiently small $\epsilon$  $g_\sigma(z)$ has the desired properties. Consider the case where $\vec\alpha=\vec 0$, i.e. $f_{\vec\alpha,\,a,\,b}(z)=a+b$. Then the function  $\tilde g_\sigma(z)=a+b+\epsilon z$ is surjective for any $\epsilon>0$. Hence it has the desired properties also in this case. Thus, $\mathcal{S}$ is dense in $E_{\rho,p}$.
\end{proof}

The following simple claim was used in the proof of Proposition \ref{prgensurj} and will be used in what follows:

\begin{lemma}\label{lmfrakCc}
Let $\mathfrak{C}$ be the set of all constant functions. Then the set $\mathfrak{C}^c=E_{\rho,p}\setminus\mathfrak{C}$ is open and dense in $E_{\rho,p}$.	
\end{lemma}
\begin{proof}
In order to prove that $\mathfrak{C}^c$ is open in $E_{\rho,p}$, we should show that the set $\mathfrak{C}$ is closed there. Consider a sequence $c_n\in\mathfrak{C}$ which converges to a function $f\in E_{\rho,p}$ in $E_{\rho,p}$, i.e. for any $n\in\N$ 
\begin{equation*}
\lim_{k\rightarrow\infty}\Vert f-c_k\Vert_n=\lim_{k\rightarrow\infty}\sup_{z\in\C}\exp\big[-(p+1/n)|z|^\rho\big]|f(z)-c_k|=0.
\end{equation*}
Since each $c_k$ is constant, this relation implies easily that $f(z)\equiv const$. Hence $\mathfrak{C}$ is closed in $E_{\rho,p}$. In order to prove that $\mathfrak{C}^c$ is dense in $E_{\rho,p}$, let us take $c\in\mathfrak{C}$ and consider the function $g_\sigma(z)=c+\sigma z\,(\sigma>0)$, which belongs to $\mathfrak{C}^c$. Then since the function $\phi(z)=z$ belongs to $E_{\rho,p}$, we obtain that for any $n\in\N$ $\lim_{\sigma\rightarrow\infty}\Vert g_\sigma-c\Vert_n=0$.
\end{proof}

In the proof of Proposition \ref{prgensurj} we have used also the following claim:

\begin{proposition}\label{lmsurjezminepz}
	The function $f(z)=e^z-\epsilon z\,(\epsilon>0)$  is surjective.
\end{proposition}
\begin{proof}
	We should prove that $f(\C)=\C$.
Take $n\in\N$ and consider the region
\begin{equation*} 
G_m=\{z\in\C:\;|\Im(z)|<\pi m,\;|\Re z|<\pi m\}.
\end{equation*}
Then 
\begin{equation*}
|f(z)|\vert_{|\Im(z)|=\pi m}=|(-1)^me^{\Re(z)}-\epsilon\Re(z)-i\epsilon\mathrm{sign}(\Im(z))\pi, m|\ge\epsilon\pi m.
\end{equation*}
\begin{equation*}
|f(z)|\vert_{\Re z=\pi m,\,\Im(z)|<\pi m}\ge e^{\pi m}-\epsilon |z|\ge e^{\pi m}-\sqrt{2}\epsilon\pi m,  
\end{equation*}
\begin{equation*}
|f(z)|\vert_{\Re z=-\pi m,\,\Im(z)|<\pi m}\ge\epsilon\pi m-1.
\end{equation*}
Hence
\begin{equation}\label{estmodfbelow}
|f(z)|\vert_{z\in\partial G_m}\ge\Theta_m,
\end{equation}
where
\begin{equation*}
\Theta_m=\min\{e^{\pi m}-\sqrt{2}\epsilon\pi m,\, \epsilon\pi m-1\}.
\end{equation*}
It is clear that
\begin{equation}\label{dfThetm}
\Theta_m\rightarrow\infty\quad\mathrm{as}\quad m\rightarrow\infty.
\end{equation}
Let us notice that by Open Mapping Theorem for analytic functions (\cite{Rud}, pp. 214, 216), the set $f(G_m)$ is open in $\C$. Furthermore, since the set $\bar G_m$ is compact in $\C$, then by Lemma \ref{lmtopol},
 \begin{equation}\label{equalpartFGpfpartG}
\partial f(G_m)=f(\partial G_m).
\end{equation}
Consider the disk 
\begin{equation}\label{dfDn}
D_{\Theta_m}=\{w\in\C:\,|w|<\Theta_m\}.
\end{equation}
It is clear that $f(0)\in D_{\Theta_m}$, if $\Theta_m>|f(0)|$. 
Let us show that
\begin{equation}\label{DnsubserfGn}
D_{\Theta_m}\subseteq f(G_m),\;\;\mathrm{if}\;\:\Theta_m>|f(0)|.
\end{equation}
Assume on the contrary that (\ref{DnsubserfGn}) is not valid. Then for some $m$, for which $\Theta_m>|f(0)|$, there is a point $w_\star\in D_{\Theta_m}\setminus f(G_m)$. Consider the closed straight segment $I_{f(0),\,z_\star}$,  lying in $D_{\Theta_m}$ and joining the points $f(0)$ and $w_\star$. Denote 
\begin{equation*}
\tilde t=\sup\{t\in[0,1]:\,f(0)+t(w_\star-f(0))\in f(G_m)\}.
\end{equation*}
Then it is clear that the point $\tilde w=f(0)+\tilde t(w_\star-f(0))$ belongs to the set
$\partial f(G_m)\cap I_{f(0),\,z_\star}$. This circumstance contradicts  (\ref{estmodfbelow}), (\ref{equalpartFGpfpartG}) and (\ref{dfDn}). Then inclusion (\ref{DnsubserfGn}) is valid, which together with \eqref{dfThetm} implies that $f(\C)=\C$.
\end{proof} 

In the proof of Proposition \ref{lmsurjezminepz} we have used the following claim:

\begin{lemma}\label{lmtopol}
Let $\mathcal{T}_k\;(k=1,2)$  be topological spaces such that $\mathcal{T}_2$ is Hausdorff, $U\subseteq\mathcal{T}_1$ be an open subset such that its closure $\bar U$ is compact  and $f:\,U\rightarrow\mathcal{T}_2$ be a continuous function such that $f(U)$ is open. Then the equality
\begin{equation}\label{eqpartials}
f(\partial U)=\partial f(U)
\end{equation}
is valid.
\end{lemma}
\begin{proof}
It is known that $f(\bar U)\subseteq\overline{f(U)}$ (\cite{Eng}, Chapt. 1, \S 4, Theorem 1, claim (v)). On the other hand since $\bar U$ is compact, then the set $f(\bar U)$ is compact in $\mathcal{T}_2$, hence it is closed there ((\cite{Eng}, Chapt. 3, \S 1, Theorems 6 and 8) . Hence $\overline{f(U)}\subseteq f(\bar U)$. Thus, the equality 
\begin{equation}\label{thusequal}
f(\bar U)=\overline{f(U)} 
\end{equation}
is valid. On the other hand, since $U$ and $f(U)$ are open in the corresponding spaces, then $\bar U=U\cup\partial U$ and $\overline{f(U)}=f(U)\cup\partial f(U)$. Then, taking into account that an open set and its boundary are disjoint, we obtain from \eqref{thusequal} the desired equality \eqref{eqpartials}.
\end{proof}

\section{Genericity of functions from $E_{\rho,p}$ with infinite number of isolated critical points}\label{sect:geninfcritval}
\setcounter{equation}{0}

The main claim of this subsection is following:

\begin{proposition}\label{prgenerinfcrit}
The set $\mathcal{N}$ of all functions from $E_{\rho,p}$ having an infinite number of isolated critical points, is $G_\delta$ and dense in $E_{\rho,p}$.
\end{proposition}
\begin{proof}
In order to prove that $\mathcal{N}$ is $G_\delta$ in $E_{\rho,p}$, we should show that the set $\mathcal{N}^c=E_{\rho,p}\setminus\mathcal{N}$ is a countable union of closed sets. By Proposition \ref{lmcndfinzer}, a non-constant function $f\in E_{\rho,p}$ belongs to t $\mathcal{N}^c$ if and only if it has the form:
\begin{equation}\label{rprfinzer}
f(z)=\int_0^zP_m(\zeta)\exp\big(Q_\rho(\zeta)\big)\,\mathrm{d}\zeta+a,
\end{equation}
where $a\in\C$ and 
\begin{equation}\label{dfpolyn1}
P_m(z)=\sum_{k=0}^mb_kz^k,\quad Q_\rho(z)=\sum_{k=1}^\rho\alpha_kz^k
\end{equation}
are polynomials  such that $\vec b=(b_0,b_1,\dots,b_m)\in\C^{m+1}$ and
\begin{equation}\label{cndforvecalph}
\vec\alpha=(\alpha_1,\alpha_2,\dots,\alpha_\rho)\in\C^{\rho-1}\times\bar D_p. 
\end{equation}
We see that $f(z)$ is not constant if and only if $\vec b\neq\vec 0$.
 
Consider the mapping $H_m$, defined by 
\begin{equation}\label{dfmapH}
H_m:\;\big(\C^{m+1}\setminus\{0\}\big)\times\big(\C^{\rho-1}\times\bar D_p\big)\times\C
\rightarrow E_{\rho,p},\quad H_m(\vec b,\,\vec\alpha,\,a)=f,
\end{equation}
where $f(z)$ has the form \eqref{rprfinzer}, and the set
\begin{equation*}
 \mathcal{C}_{m,\,l}=,\big(\bar D_l^{m+1}\setminus\bar D_{\delta_l}\big)\times\big(\bar D_l^{\rho-1}\times\bar D_p\big)\times\bar D_l,
\end{equation*}
 where $\delta_l=1/(l+1)$ . Notice that by Proposition \ref{lmcontHba} the mappings $H_m$ are continuous. Therefore the sets $\mathcal{K}_{m,\,l}=H_m(\mathcal{C}_{m,\,l})$ are compact in the space $E_{\rho,p}$ as images of compact sets, hence they are closed in $E_{\rho,p}$. Thus, we obtain the following representation for the set $\mathcal{N}^c$:
\begin{equation*}
\mathcal{N}^c=\big(\bigcup_{m=0}^\infty\bigcup_{l=1}^\infty\mathcal{K}_{m,\,l}\big)\cup\mathfrak{C}. 
\end{equation*}
Recall that $\mathfrak{C}$ is the set of all constant functions, which by Lemma \ref{lmfrakCc} is closed in $E_{\rho,p}$. These circumstances mean that the set $\mathcal{N}$ is $G_\delta$ in $E_{\rho,p}$.

Let us show that $\mathcal{N}$ is dense in $E_{\rho,p}$. Assume that $f\in\mathcal{N}^c$. If $\vec\alpha\neq \vec 0$, we see from \eqref{rprfinzer} that $f^\prime (z)$ is not a polynomial. Then by Great Picard Theorem, the function $f^\prime(z)$ takes every complex value, with at most one exception, infinitely
many times (\cite{Gir}). If this exceptional value is not zero, this function has infinite number of isolated zeros, which are critical points of $f$. Otherwise for any $\epsilon>0$ the function $f^\prime(z)+\epsilon$ has infinite number of zeros, which are critical points of the function $g(z)=f(z)+\epsilon z$ belonging to $E_{\rho,p}$. It is clear that for any $n\in\N$ and $\sigma>0$ there is $\epsilon>0$ such that $\Vert f-g\Vert_n<\sigma$. Now consider the case where $\vec\alpha=\vec 0$, i.e. $f(z)$ is a polynomial. Then for $\beta\in\bar D_p\setminus\{0\}$ and $\eta>0$ the function $h(z)=f(z)+\eta\exp(\beta z)$ belongs to $E_{\rho,p}$ and its derivative $h^\prime(z)$ is not a polynomial. It is obvious that for any $n\in\N$ and $\sigma>0$ there is $\eta>0$ such that $\Vert f-h\Vert_n<\sigma/2$. Applying again Great Picard Theorem to the function $h^\prime(z)$, we can find $\epsilon\ge 0$ such that the function $g(z)=h(z)+\epsilon z$ belongs to $E_{\rho,p}$, has an infinite number of isolated critical points and $\Vert g-h\Vert_n<\sigma/2$. Hence $\Vert f-g\Vert_n<\sigma$. Thus, $\mathcal{N}$ is dense in $E_{\rho,p}$.
\end{proof}

\section{Proof of main results}\label{sect:proofmainres}\label{sect:proofmainres}
\setcounter{equation}{0}

\subsection{Proof of Theorem \ref{mainthpol}}
\begin{proof}
	(i) In view of Definition \ref{defdisccritval} and Remark \ref{remcritval}
	, condition (\ref{maincond}) implies that all the
	critical values $P_m(y_1),P_m(y_2),\dots, P_m(y_{m-1})$ of $P_m$ are pairwise different. Hence, in particular, all the critical points
	$\,y_1,y_2,\dots,y_{m-1}$ of $P_m$  are pairwise different, hence they are
	non-degenerate. 
	
	(ii) It is clear that the Riemann surface of $P_m^{-1}$ is connected. Then claim (i) and Lemma 1 of \cite{Z} imply that the the group ${\mathrm Mon}(\vec a)$ is transitive and
	generated by transposition. Then by Lemma 2 of \cite{Z}, $\mathrm {Mon}(\vec a)=\mathrm{S}(m)$.
\end{proof}

\subsection{Proof of Theorem \ref{thgenfrakG}}
\begin{proof}
Recall that $\mathfrak{C}$ is the set of all constant functions and by Lemma \ref{lmfrakCc} the set $\mathfrak{C}^c=E_{\rho,p}\setminus\mathfrak{C}$ is open and dense in $E_{\rho,p}$. Furthermore, and all the functions from $\mathfrak{C}^c$ have isolated critical points. Let us take $f\in\mathfrak{C}^c$ and  $l\in\N$. Let ${\mathcal Cr}_f=\{z_{k,\,f}
\}_{k=1}^{K_f}$ be the set of all the critical points of $f$, numbered in non-decreasing order of their modules ($K_f$ can be both finite and infinite). Let $\{\mathcal{P}_{f,\,m}\}_{m=1}^{M_f}$ be the partition of ${\mathcal Cr}_f$, generated by the following equivalence relation: $z_{\mu,\,f}\sim z_{\nu,\,f}$ if and only if $|z_{\mu,\,f}|= |z_{\nu,\,f}|$. Let $r_{f,\,m}$ be the radius of the circle $C_{f,\,m}$ with the center at $0$, such that $\mathcal{P}_{f,\,m}\subset C_{f,\,m}$. Then the sequence $\{r_{f,\,m}\}_{m=1}^{M_f}$
is increasing. Consider the quantities:
\begin{eqnarray*}
s_{f,\,m}=\left\{\begin{array}{ll}
\frac{r_{f,\,m}+r_{f,\,m+1}}{2},&\mathrm{if}\quad m\le M_f-1,\\
r_{f,\,m}+1,&\mathrm{if}\quad m= M_f
\end{array}
\right.
\end{eqnarray*}
and
\begin{eqnarray*}
R_{f,\,l}=\left\{\begin{array}{ll}
l,&\mathrm{if}\quad \sup_{1\le m\le M_f}s_m\le l,\\
\min_{m:\,s_m\ge l}\,s_m&\mathrm{otherwise}.
\end{array}
\right.
\end{eqnarray*}
It is clear that  each circle $\tilde C_{f,\,l}$, whose  center is $0$ and radius  $R_{f,\,l}$, contains no critical points of $f$. Let $D_{f,\,l}$ be an open disk, whose boundary is $\tilde C_{f,\,l}$ and $D_L$ is an open disk with the center $0$ and radius $I$. Since $R_{f,\,l}\ge l$, $D_l\subseteq D_{f,\,l}$ for any $f\in\mathfrak{C}^c$, Denote by $\mathfrak{V}_{l}$ (by $\mathfrak{V}$) the set of all functions from $\mathfrak{C}^c$ having only non-degenerate critical points in $D_l$ (respectively, in $\C$)  and taking at them pairwise different values. It is clear that 
\begin{equation}\label{intersecVl}
\mathfrak{V}=\bigcap_{l\in\N}\mathfrak{V}_{l}.
\end{equation}
By Proposition \ref{prlocgener}, for any $f\in\mathfrak{C}^c$ there is its open connected neighborhood $U(f)\subseteq\mathfrak{C}^c$ and a principal analytic variety $\mathcal{V}_f\subset U(f)$ (Definition \ref{dfprincvar}
) such that the set of all functions $w\in U(f)$ having only non-degenerate critical points in $D_{f,\,l}$ and taking at them pairwise different values, coincides with the set $U(f)\setminus\mathcal{V}_f$, which is open in $U(f)$ (hence in $E_{\rho,p}$). Furthermore, since the variety $\mathcal{V}_f$ is an analytic proper subset of $U(f)$, it is nowhere dense in $U(f)$ by Proposition \ref{pruniq} (which is a generalization of uniqueness theorem for holomorphic functions).
Then the set $U(f)\setminus\mathcal{V}_f$ is dense in $U(f)$. Notice that since  $D_l\subseteq D_{f,\,l}$ for any $l\in\N$ and $f\in\mathfrak{C}^c$, then $U(f)\setminus\mathcal{V}_f\subseteq\mathfrak{V}_{l}$. For each $l\in\N$ consider the set 
\begin{equation*}
\tilde{\mathfrak V}_l=\bigcup_{f\in\mathfrak{C}^c}U(f)\setminus\mathcal{V}_f.
\end{equation*}
It is clear that $\tilde{\mathfrak V}_l\subseteq\mathfrak{V}_{l}$ and $\tilde{\mathfrak V}_l$ is open and dense in $\mathfrak{C}^c$. On the other hand, since  the set $\mathfrak{C}^c$ is open and dense in $E_{\rho,p}$, the set $\tilde{\mathfrak V}_l$ is open and dense in $E_{\rho,p}$.  Then the set 
\begin{equation*}
\tilde{\mathfrak V}=\bigcap_{l\in\N}\tilde{\mathfrak V}_l
\end{equation*}
is $G_\delta$ in $E_{\rho,p}$, hence by Baire Category Theorem it is dense in $E_{\rho,p}$. Furthermore, in view of \eqref{intersecVl}, $\tilde{\mathfrak V}\subseteq{\mathfrak V}$. Recall that by Propositions \ref{prgensurj} and \ref{prgenerinfcrit}, the set $\mathcal{S}$ of all surjective functions from $E_{\rho,p}$ and the set $\mathcal{N}$ of all functions from $E_{\rho,p}$ having infinite number of isolated critical points are $G_\delta$ and dense in $E_{\rho,p}$. Then the set $\mathfrak{T}=\mathfrak{V}\cap\mathcal{S}\cap\mathcal{N}$ contains the set $\tilde{\mathfrak V}\cap\mathcal{S}\cap\mathcal{N}$, which is $G_\delta$ and dense in $E_{\rho,p}$.
\end{proof}

\subsection{Proof of Theorem \ref{prmongrfottyp}}

\begin{proof}
The way of the proof will be similar to one in \cite{Z}, \S 4, Lemma 1. Let $f(w)$ be a function from $\mathfrak{T}$, where the case of algebraic variety was considered. Then it is surjective, has infinite number of non-degenerate critical points $w_k\,(k\in\N)$ and its critical values $z_k=f(w_k)$ are pairwise different. Hence in particular, $f(w)$ is not a polynomial, therefore by Great Picard Theorem, for any
$z\in\C$ the set $f^{-1}(z)$ is countable. Consider the sets
\begin{equation*}
\Phi=\{(z,\,w)\in\C^2:\,z=f(w) \},\quad
\mathcal{F}=\bigcup_{k\in\N}\{(z_k,\,w_k)\},\quad \mathfrak{G}=\Phi\setminus\mathcal{F}.
\end{equation*}
Then $f:\,\Phi\rightarrow\C$ is a ramified covering of the plain $\C$ by $\Phi$ (\cite{Sp}, Chapt 4) and $(z_k,\,w_k)\,(k\in\N)$ are its branching points. It is clear that the set $\mathfrak{G}$ can be identified with Riemann surface of the multi-valued function $f^{-1}$, which consists of countable number of sheets. Let is numerate them by naturals $\N$. Consider the set of all
loops $\gamma$ in $f(\mathfrak{G})$ with beginning and end at a point $a\in f(\mathfrak{G})$ and
consider the corresponding monodromy group $\mathrm{Mon}(\Phi)$. It
consists of the bijective mappings
$\Delta_\gamma:\;\N\rightarrow\N$, such that each of them
corresponds
to a class of equivalent loops (with respect
to deformations in $f(\mathfrak{G})$).  Observe that each $(z_k, w_k)$  is the branching point of the second
order, because $f^{\prime\prime}(w_k)\neq 0$. Then since the points $z_k$ are pairwise different, the group
$\mathrm{Mon}(\Phi)$ is generated by transpositions. On the other
hand, $\mathrm{Mon}(\Phi)$ acts on $\N$ transitively, because
the Riemann surface $\mathfrak{G}$
is connected. Hence by Lemma \ref{lmfinsymgr},  $\mathrm{Mon}(\Phi)=\mathrm{FSym}(\N)$.
\end{proof}

In the proof of Theorem \ref{prmongrfottyp} we have used the following claim: 

\begin{lemma}\label{lmfinsymgr}
If a transitive group $G$ of permutations on $\N$ is generated by transpositions, it coincides with $\mathrm{FSym}(\N)$.
\end{lemma}
\begin{proof}
Since $G$ is generated by transpositions, it is a subgroup of $\mathrm{FSym}(\N)$ and by claim (b) of Theorem 5.2 from \cite{Hall}, $G=\bigoplus_{i\in I}\mathrm{FSym}(\Omega_i)$, where $\Omega_i\,(i\in I)$ are all the distinct orbits of $G$. But since $G$ is transitive, the unique orbit of $G$ is the whole $\N$. Thus, $G=\mathrm{FSym}(\N)$.
\end{proof}

\subsection{Proof of Theorem \ref{prcon1}}

\begin{proof} Let  $\mathrm{Alt}(\N)$ be the subgroup of the group $\mathrm{FSym}(\N)$
	consisting of finitary permutations which are compositions of even
	number of transpositions (the alternating subgroup). Since $\mathrm{Alt}(\N)$
	is a normal  subgroup of $\mathrm{FSym}(\N)$ with the two-element quotient
	group $\mathrm{FSym}(\N)/\mathrm{Alt}(\N)$, it is enough to show that $\mathrm{Alt}(\N)$ is interspersed-finite unsolvable.  To
	this end assume on the contrary that it is solvable in this sense, i.e. by Definition \ref{defofabfinsolv} we have the sequence \eqref{seqofgr} of subgroups $G_k$ of the group $G=\mathrm{Alt}(\N)$, 
	such that for each $k\in\{0,\,1,\,2,\,\dots,\,r-1\}$ $G_{k+1}$ is a normal subgroup of $G_k$ and the quotient group $G_k/G_{k+1}$ is either abelian or finite. But it is known that the group $\mathrm{Alt}(\N)$ is simple (\cite{Sc}, p. 295), i.e. it cannot contain nontrivial normal subgroups. Hence it must be only the case where $G_1=\{e\}$ (i.e,, r=1). But since $\mathrm{Alt}(\N)=\mathrm{Alt}(\N)/G_1$ is infinite and non-abelian, we obtain  contradiction to the assumption. 
	\end{proof}

\subsection{Proof of Corollary \ref{corrol}}

\begin{proof}
Assume on the contrary that the function $f^{-1}$ is representable in the manner indicated in the formulation of this  corollary. Then in view of \cite{Kh}, p. 720 (Results on generalized quadratures) the monodromy group of its Riemann surface is interspersed-finite solvable. This contradicts  Theorems \ref{prmongrfottyp} and \ref{prcon1}. 	
\end{proof}

\section{Some problems}\label{sect: pronlems}
\setcounter{equation}{0}

\subsection{A problem connected with Theorem \ref{mainthpol}}\label{subsect: pronlem1}

In \cite{Katz} a stratification of the discriminant variety $D_2$ 
\begin{equation*}
D_2\supset D_3\supset\dots\supset D_m
\end{equation*}
for polynomials $P_m(z)=z^m+a_{m-1}z^{m-1}+\dots+a_0$ is considered, where $D_k$ consists of those points of the coefficient space 
$\C^m$, which correspond to polynomials having at least one root with multiplicity $\ge k$. It was established in \cite{Katz} that $D_{k-1}$ is comprised of the affine subspaces tangent to $D_{k}$,
and the number of such subspaces which hit a given point $P_m\in D_{k-1}$ is entirely
determined by the multiplicities of the $P_m(z)$-roots. We can consider the analogous stratifications of the discriminant variety $V^\prime$ of the derivative $P_m^\prime(z)$ and of $CV$-discriminant variety $V_{cvd}$, which is the discriminant variety of the polynomial 
\begin{equation*}
U_{m-1}^{\vec
	a}(y)=\mathrm {det}(P_m(C(\vec a))-yI)\quad (\vec a=(a_0, a_1,\dots,a_{m-1}))
\end{equation*}
(Definition \ref{defdisccritval} and Remark \ref{remcritval}). Since each stratum of $V_{cvd}$ is connected with the number of cycles taking part in generators of monodromy group $\mathrm{Mon}(\vec a)$ of Riemann surface of $P_m^{-1}$ and each stratum of $V^\prime$ is connected with orders of these cycles, the following problem appears:

\begin{problem}\label{prob1}
To what extent the stratifications of $V_{cvd}$ and $V^\prime$ and affine tangent subspaces to their strata define the monodromy group $\mathrm{Mon}(\vec a)$ ?  Is it possible to formulate a condition for isomorphism of $\mathrm{Mon}(\vec a_1)$ and $\mathrm{Mon}(\vec a_2)$ in these terms ?
\end{problem}

\subsection{A problem connected with Theorem \ref{thgenfrakG}}\label{subsect: pronlem2}
	
Consider on the complex vector space $E$ of all entire functions the topology of uniform convergence on compacts of the complex plane $\C$. Then the set  of all the polynomials $P_m(z)$ of degree $m\ge 1$ is dense in $E$. On the other hand, each such a polynomial does not belong to the set $\mathfrak{T}$ of typical functions, since it has only finite number of critical points. Therefore the set $E\setminus\mathfrak{T}$ is dense in $E$. Thus, in this situation Theorem \ref{thgenfrakG} is not true.
Hence the following problem appears:

\begin{problem}\label{prob2}
Is there a vector space $\tilde E$ of entire functions with a Frechet topology generated by a sequence of weighted norms, which is a proper extension of all the spaces $E_{\rho,p}$, such that Theorem \ref{thgenfrakG} is valid in it ? To what extent can the spaces $E_{\rho,p}$ be expanded in this way ?
\end{problem}

\appendix

\section{Auxiliary claims}\label{sect:appendix}
\setcounter{equation}{0}

\subsection{Splitting of zeros \and interpolation}\label{subsect:splitinterp}

In the proof of Proposition \ref{prdifcritval} we have used the following claims:

\begin{proposition}\label{lmsplitrots}
	Let $w_0(z)$ be an entire function from a weighted space $H(v)$
	Assume that $w_0(z)$  has in the disk $D_R$ $m$ zeros $z_1^0,\,z_2^0,\,\dots,\,z_m^0$ (counting their multiplicities) and it does not vanish on $\partial D_R$. Then for any $\epsilon>0$ there is a function $w_\epsilon\in H(v)$ such that $w_\epsilon-w_0$ belongs to the space  $H(\tilde v)\subset H(v)$ with the weight 
	\begin{equation}\label{dfweighv}
	\tilde v(r)=(1+r^2)v(r),
	\end{equation}
	$\Vert w_\epsilon-w_0\Vert_{\tilde v}<\epsilon$ 
	and $w_\epsilon(z)$  has in $D_R$ $n$ pairwise different zeros and it does not vanish on $\partial D_R$. 
\end{proposition}
\begin{proof}
	Let $z_1^0,\,z_2^0,\,\dots,\,z_p^0$ be all pairwise different zeros of $w_0(z)$. Denote by $m_k$ the multiplicity of $z_k^0$ and assume that all the multiple zeros (i.e. with $m_k>1$) are $\{z_k^0 \}_{k=1}^l\,(1\le l\le p)$. Let us show that by a $H(\tilde v)$-small perturbation it is possible to split the zero $z_1$ into $m_1$ pairwise different zeros, not varying the rest of them, without appearance of new zeros in the closed disk $\bar D_R$.
	Consider the representation 
	\begin{equation}\label{reprw0z}
	w_0(z)=P_m(z)\Phi(z),
	\end{equation}
	where 
	\begin{equation}\label{dfPnmult}
	P_m(z)=\prod_{k=1}^p(z-z_k^0)^{m_k}
	\end{equation}
	and $\Phi(z)$ is an entire function.
	We construct the desired perturbation of $w_0(z)$in the form $w_\delta^{(1)}(z)=w_0(z)-\delta\Psi(z)$ with $\delta>0$ and 
	\begin{equation}\label{dfPsiz}
	\Psi(z)=\prod_{k=2}^p(z-z_k^0)^{m_k}\Phi(z). 
	\end{equation} 
	Then we have:
	\begin{equation}\label{anotherdfwdel}
	w_\delta^{(1)}(z)=[(z-z_0)^{m_1}-\delta]\Psi(z).
	\end{equation}
	Notice that in view of \eqref{reprw0z}, \eqref{dfPnmult} and \eqref{dfPsiz}, for some $C>0$, $N>0$
	\begin{equation}\label{estforPsi}
	|\Psi(z)|\le C\frac{|w_0(z)|}{|z|^{m_1}},\quad\mathrm{if}\quad |z|\ge N.
	\end{equation}
	Since $w_0\in H(v)$ and $m_1>1$, the last estimate implies that $\Psi\in H(\tilde v)$ with the weight $\tilde v$, defined by \eqref{dfweighv}. Since $H(\tilde v)\subset H(v)$, $w_\delta^{(1)}\in H(v_)$. Let us take an arbitrary $\epsilon>0$.
	We  see from \eqref{anotherdfwdel} that  $m_1$ pairwise different numbers $z_{\nu,\i}=\delta^{1/m_1}\exp\big(\frac{i\,2\pi}{m_1}\nu\big)\,(\nu=0,1,2,\dots,m_1-1)$ are zeros of the function  $w_\delta^{(1)}(z)$ and if $\delta^{1/m_1}<\min_{2\le k\le p}|z_k^0-z_1^0|$, the numbers  $z_k^0\,(k=2,\dots, p)$ are its zeros  with the  multiplicities $m_k$. Furthermore, for $\delta$ small enough $\Vert w_\delta^{(1)}-w_0\Vert_{\tilde v}<\epsilon/l$ and, by Rouchet's Theorem, the function $w_\delta^{(1)}(z)$ has in the closed disk $\bar D_R$ only the zeros, mentioned above. Afterwards the proof is finished by induction. 
\end{proof}

\begin{proposition}\label{lmintprob}
	Let $z_k\,(k=1,2,\dots,m)$ be pairwise different numbers and $M_{2m}
	(z,j)$ be polynomials of the form
	\begin{eqnarray}\label{dffoundpol}
	M_{2m}(z,j)=[(z-z_j)^2+1]\frac{\prod_{k\in\{1,2,\dots,m\}\setminus\{j\}}(x-x_k)^2}{\prod_{k\in\{1,2,\dots,m\}\setminus\{j\}}(x_j-x_k)^2}
	\end{eqnarray}
	Then the polynomial
	\begin{equation}\label{dfintpol}
	P_{2m}(z)=\sum_{j=1}^ny_jM_{2m}(z,j)
	\end{equation}
	is a solution of the interpolation problem
	\begin{equation*}
	P_{2m}(z_j)=y_j,\quad P_{2m}^\prime(z_j)=0\quad (j=1,2,\dots,m). 
	\end{equation*}
\end{proposition}
\begin{proof}
	We see from \eqref{dffoundpol} that $M_{2m}(z_k,j)=\delta_{j,k}$  and $\frac{\mathrm{d}}{\mathrm{d}z}M_{2m}(z_k,j)=0$ $(j,k\in\{1,2,\dots,m\})$. This fact and \eqref{dfintpol} imply the desired claim. 
\end{proof}

\subsection{Properties of some mappings}\label{subsect:propsomemap}

In the proof of Proposition \ref{prgensurj} we have used the following claims:

\begin{proposition}\label{lmpropfamexp}
	
	(i) Any non-surjective function from $E_{\rho,p}$ has the form \eqref{dffalpha};
	\vskip2mm

\noindent (ii)  the mapping $\Phi:\;(\C^{\rho-1}\times\bar D_p)\times\C\times\C\rightarrow E_{\rho,p}$, defined by $\Phi(\alpha,\,a,\,b):=f_{\vec\alpha,\,a,\,b }$, is continuous.
	\vskip2mm
	
\end{proposition}
\begin{proof}
	(i) A function $f\in E_{\rho,p}$ is not surjective if and only if there is a number $a\in C$, which does not belong to $f(\C)$, i.e. $f(z)-a\neq 0$ for any $z\in\C$. Hence the function $\phi(z)=\frac{f^\prime(z)}{f(z)-a}$ is entire. This means that $f(z)$ has the form
	\begin{equation}\label{reprsurjfunc}	
	f(z)=c\exp(\psi(z))+a,
	\end{equation}
	where $c,a\in\C$ and $\psi(z)=\int_0^z\phi(s)\,\mathrm{d}s$ is an entire function. It is easy to see that function of the form \eqref{reprsurjfunc} is non-surjective. On the other hand, by definition \eqref{dfexptyp} of the space $E_{\rho,p}$, for any $\epsilon>0$ there is $R>0$ such that
	\begin{equation}\label{usedfexptyp}
	\frac{\ln|\exp(\psi(z))|}{|z|^\rho}< p+\epsilon \quad\forall\; z\in\C\setminus\bar D_R.
	\end{equation}
	Hence since the function $\exp(\psi(z))$ has no zeros, then by Hadamard's Factorization Theorem for entire functions of finite order, the function $\psi(z)$ is a polynomial of the form $\psi(z)=\sum_{k=0}^\rho \alpha_kz^k$ (\cite{Mar}, II.10, p. 289).  We obtain from \eqref{reprsurjfunc}  the desired representation \eqref{dffalpha} with $b=ce^{\alpha_0}$. Then inequality \eqref{usedfexptyp} is equivalent to 
\begin{equation*}
\forall\,\epsilon>0\;\exists\, R>0\;\forall\; z\in\C\setminus\bar D_R:\quad \Re\Big(\alpha_\rho\frac{z^\rho}{|z|^\rho}\Big)+ \Re\Big(\sum_{k=1}^{\rho-1} \alpha_k\frac{z^k}{|z|^\rho}\Big)<p+\epsilon.
\end{equation*}
Taking into account that $\epsilon>0$ is arbitrary and representing $z=|z|e^{i\varphi}$ and $\alpha_\rho=|\alpha_\rho|e^{i\theta}\;(\varphi,\theta\in[0,\,2\pi))$, we obtain that the last inequality is equivalent to:
$|a_\rho|\cos(\rho\varphi+\theta)\le p$ for any $\varphi\in[0,\,2\pi)$. It is clear that the last fact is equvalent to $|a_\rho|\le p$. Claim (i) is proven.

	  (ii) We should show that $\Phi$ is continuous with respect to each norm $\Vert\cdot\Vert_n$ on $E_{\rho,p}$, defined by \eqref{dfnormn}. To this end consider the functions
	\begin{equation*}
	F_n(\vec\alpha,\,a,\,b,\,z)=\exp(-(p+1/n)|z|^\rho)\Big(b\exp\big(\sum_{k=1}^\rho \alpha_kz^k\big)+a\Big). 
	\end{equation*}
	Let us take $R>0$, $\epsilon>0$ and $\vec\alpha_0=(a_1^0, a_2^0,\dots, a_\rho^0)$. Then
	\begin{eqnarray}\label{estdiffer}
	&&\hskip-10mm\Vert \Phi(\vec\alpha,\,a,\,b)-\Phi(\vec\alpha_0,\,a_0,\,b_0)\Vert_n=\Vert f_{\vec\alpha,\,a,\,b}-f_{\vec\alpha_0,\,a_0,\,b_0}\Vert_n\le\nonumber\\
	&&\hskip-10mm\max\big\{\sup_{z\in D_R}|F_n(\vec\alpha,\,a,\,b,\,z)-F_n(\vec\alpha_0,\,a_0,\,b_0\,z)| ,\nonumber\\
	&&\hskip-10mm\sup_{z\in\C\setminus D_R}\Big(|b|\exp\big(\sum_{k=1}^\rho|\alpha_k||z|^k-(p+1/n)|z|^\rho\big)+\\
	&&\hskip-10mm|b_0|\exp\big(\sum_{k=1}^\rho|\alpha_k^0||z|^k-(p+1/n))|z|^\rho)\big)+\nonumber\\
	&&\hskip-10mm\big(|a|+|a_0|\big)\exp\big(-(p+1/n)|z|^\rho\big)\Big)\big\}.\nonumber
	\end{eqnarray}
	Assume that $\vec\alpha,\vec\alpha_0\in\bar D_r^{\rho-1} \times\bar D_p\;(r>0)$, $a,a_0\in\bar D_r$ and $b,b_0\in\bar D_r$. Let us choose $R>0$ such that
	\begin{eqnarray}\label{estoftail}
	&&\sup_{z\in\C\setminus D_R}\Big(|b|\exp\big(\sum_{k=1}^\rho|\alpha_k||z|^k-(p+1/n)|z|^\rho\big)+\nonumber\\
	&&|b_0|\exp\big(\sum_{k=1}^\rho|\alpha_k^0||z|^k-(p+1/n))|z|^\rho\big)+\nonumber\\
	&&\big(|a|+|a_0|\big)\exp\big(-(p+1/n)|z|^\rho\big)\Big)\le\\
	&&2r\exp\big(r\sum_{k=1}^{\rho-1}R^k-(1/n)R^\rho\big)+2r\exp\big(-(p+1/n)R^\rho\big)<\epsilon.\nonumber
	\end{eqnarray} 
	On the other hand, since the function $F_n(\vec\alpha,\,a,\,b,\,z)$
	is continuous on the compact $(\bar D_r^{\rho-1} \times\bar D_p)\times\bar D_r\times\bar D_r\times\bar D_R$, it is uniformly continuous there. Hence there is $\delta>0$ such that if $\vec{a}\alpha,\vec\alpha_0\in \bar D_r^{\rho-1} \times\bar D_p$, $a,a_0\in\bar D_r$, $b,b_0\in\bar D_r$, $|\vec\alpha-\vec\alpha_0|<\delta$, $|a-a_0|<\delta$ and $|b-b_0|<\delta$, then
	\begin{equation*}
	\sup_{z\in D_R}|F_n(\vec\alpha,\,a,\,b,\,z)-F_n(\vec\alpha_0,\,a_0,\,b_0,\,z)|<\epsilon.
	\end{equation*}
	 The last circumstance and estimates \eqref{estdiffer}, \eqref{estoftail} imply that 
	 \begin{equation*}
	 \Vert\Phi(\vec\alpha,\,a,\,b)-\Phi(\vec\alpha_0,\,a_0,\,b_0)\Vert_n<\epsilon.
	 \end{equation*}
	 Since $r>0$ is arbitrary, then the mapping $\Phi$ is continuous on the whole its domain of definition. Claim (ii) is proven.
	\end{proof}

In the proof of Proposition \ref{prgenerinfcrit} we have used the following claims:

\begin{proposition}\label{lmcndfinzer}
	A function $f\in E_{\rho,p}$ has at most finite number of critical points if and only if it has the form \eqref{rprfinzer}-\eqref{dfpolyn1}, where $\vec\alpha$ satisfies condition \eqref{cndforvecalph}.
\end{proposition}
\begin{proof}
	Suppose that the function $f\in E_{\rho,p}$ has at most finite number of critical points, i.e. its derivative $f^\prime(z)$ has $m$ zeros (counting their multiplicities). It is clear that this fact is equivalent to the representation $f^\prime(z)=R_m(z)\phi(z)$, where $R_m(z)$ is a polynomial of degree $m$ and $\phi(z)$ is an entire function without zeros, i.e. the function $\frac{\phi^\prime(z)}{\phi(z)}$ is entire. Therefore, using arguments from the proof of claim (i) of Proposition \ref{lmpropfamexp}, we have the representation
	\begin{equation*}
	f^\prime(z)=R_m(z)\exp(\varphi(z)),
	\end{equation*}
	where $\varphi(z)$ is an entire function. On the other hand, since $f\in E_{\rho,p}$, then by claim (i) of Lemma \ref{lmholofwprim}, $f^\prime\in E_{\rho,p}$, i.e.  
	\begin{equation*}
	\limsup_{|z|\rightarrow\infty}\frac{\ln|R_m(z)\exp(\varphi(z))|}{|z|^\rho}\le p.
	\end{equation*}
	It is clear that this inequality is equivalent to
	\begin{equation*}
	\limsup_{|z|\rightarrow\infty}\frac{\ln|\exp(\varphi(z))|}{|z|^\rho}\le p.
	\end{equation*}
	Using again arguments from the proof of claim (i) of Proposition \ref{lmpropfamexp}, we obtain 
	  that the last inequality is equivalent to the representation $\varphi(z)=\sum_{k=0}^\rho\alpha_kz^k$ with the vector $\vec\alpha=(\alpha_1,\alpha_2,\dots,\alpha_\rho)$ satisfying condition \eqref{cndforvecalph}. Thus, after integration we obtain representation \eqref{rprfinzer}-\eqref{dfpolyn1}, where
	 \begin{equation*}
	 P_m(z)=e^{\alpha_0} R_m(z).\qedhere 
	 \end{equation*} 
	   . \end{proof}

\begin{proposition}\label{lmcontHba}
	The mapping $H_m$, defined by \eqref{dfmapH} and \eqref{rprfinzer}, is continuous.
	\end{proposition}
\begin{proof} We should show that the mapping $H_m$ is continuous with respect to each norm $\Vert\cdot \Vert_n$ on $E_{\rho,p}$, defined by \eqref{dfnormn}. Let us take 
	\begin{eqnarray*}
	&&\vec b=(b_{0},b_{1},\dots,b_{m}),\; \vec b_0=(b_{0}^0,b_{1}^0,\dots,b_{m}^0)\in\C^{n+1},\nonumber\\
	&&\vec\alpha=(\alpha_{1},\alpha_{2},\dots,\alpha_{\rho}),\; \vec \alpha_0=(\alpha_{1}^0,\alpha_{2}^0,\dots,\alpha_{\rho}^0)\in\C^{\rho},\nonumber\\
	&&a,a_0\in\C,
	\end{eqnarray*}
	$ R>0$ and  $\epsilon>0$. Then
	\begin{eqnarray}\label{estnormdifH}
	&&\hskip-12mm\Vert H_m(\vec b,\,\vec\alpha,\, a)-H_m(\vec b_0,\,\vec\alpha_0,\, a_0)\Vert_n\le\nonumber\\
	&&\hskip-12mm\max\{\sup_{z\in D_R}\big(\exp(-(p+1/n)|z|^\rho)|H_m(\vec b,\,\vec\alpha,\, a)(z)-
	H_m(\vec b_0,\,\vec\alpha_0,\, a_0))(z)|\big),\nonumber\\
	&&\hskip-12mm\;\sup_{z\in\C\setminus D_R}\big(\exp(-(p+1/n)|z|^\rho)(|H_m(\vec b,\,\vec\alpha,\, a)(z)|+\nonumber\\
	&&|H_m(\vec b_0,\,\vec\alpha_0,\, a_0)(z)|)\big)\}.
	\end{eqnarray} 
	Assume that 
	\begin{equation}\label{assthat}
	\vec b,\,\vec b_0\in\bar D_r^{m+1}\,(r>0),\quad \vec\alpha,\,\vec\alpha_0\in\bar D_r^{\rho-1}\times\bar D_p,\quad a,\, a_0\in \bar D_r. 
	\end{equation}
	Using \eqref{dfmapH}, \eqref{rprfinzer}, \eqref{dfpolyn1}, we have:  
	\begin{eqnarray*}
		&&\sup_{z\in\C\setminus D_R}\big(\exp\big(-(p+1/n)|z|^\rho\big)|H_m(\vec b,\,\vec\alpha,\, a)(z)|\big)\le\nonumber\\
		&&\sup_{z\in\C\setminus D_R}\Big(\exp\big(-(p+1/n)|z|^\rho\big)\big(\sum_{k=0}^mr\int_0^{|z|}t^k
		\exp\big(pt^\rho+\\
		&&r\sum_{j=1}^{\rho-1} t^j\big)\,\mathrm{d}t+a\big)\Big)
	\end{eqnarray*}
Using L'Hospital's rule, we obtain:
	\begin{eqnarray*}
	&&\hskip-18mm\lim_{|z|\rightarrow\infty}\exp\big(-(p+1/n)|z|^\rho\big)\sum_{k=0}^m\int_0^{|z|}t^k
	\exp\big(pt^\rho+r\sum_{j=1}^{\rho-1} t^j\big)\,\mathrm{d}t=\\
	&&\hskip-18mm\lim_{|z|\rightarrow\infty}\frac{\sum_{k=0}^m|z|^k
		\exp\big(p|z|^\rho+r\sum_{j=1}^{\rho-1} |z|^j\big)}{\exp\big((p+1/n)|z|^\rho\big)\rho|z|^{\rho-1}}=0.
	\end{eqnarray*}
	These circumstances mean that we can choose $R>0$ so that
	\begin{equation}\label{estoftail}
	\sup_{z\in\C\setminus D_R}\big(\exp(-(p+1/n)|z|^\rho)(|H_m(\vec b,\,\vec\alpha,\, a)(z)|\big)<\epsilon/2.
	\end{equation}
	On the other hand, in view of \eqref{dfmapH}, \eqref{rprfinzer}, \eqref{dfpolyn1}, the function
	\begin{equation*}
	\exp(-(p+1/n)|z^\rho|)H_m(\vec b,\,\vec\alpha,\, a)(z)
	\end{equation*}
	is continuous with respect to all the variables $\vec b,\,\vec\alpha,\,a,\,z$, hence it is uniformly continuous on the compact $\bar D_r^{m+1}\times(\bar D_r^{\rho-1}\times\bar D_p)\times\bar D_r\times\bar D_R$. Hence there is $\delta>0$ such that 
	\begin{equation*}
	\sup_{z\in D_R}\big(\exp(-(p+1/n)|z|^\rho)|H_m(\vec b,\,\vec\alpha,\, a)(z)-
	H_m(\vec b_0,\,\vec\alpha_0,\, a_0))(z)|\big)<\epsilon
	\end{equation*}
	as soon as the condition \eqref{assthat} is fulfilled and 
	\begin{equation*}
	\max\{|\vec b-\vec b_0|,\,|\vec\alpha-\vec\alpha_0|,\,|a-a_0|\}<\delta.
	\end{equation*}
	The last circumstance and estimates \eqref{estnormdifH}, \eqref{estoftail} imply that under conditions mentioned above $\Vert H_m(\vec b,\,\vec\alpha,\, a)-H_m(\vec b_0,\,\vec\alpha_0,\, a_0)\Vert_n<\epsilon$. Since $r>0$ is arbitrary, then the mapping $H_m$ is continuous on the whole its domain of definition. 
\end{proof}

\subsection{Holomorphic functions on complex Frechet spaces}\label{subsect:holfuncFrech}

\begin{proposition}\label{pronintpar}
Let $\Omega$ be an open subset of a complex Frechet space $E$ and  $D$ be a bounded open domain in the complex plane $\C$ such that its boundary $\partial D$ is rectifiable. Let $g(w,z)$ be a function defined  in $\Omega\times\partial D$ and continuous there, such that it is holomorphic with respect to $w$ in $\Omega$ for each fixed $z\in\partial D$. Then the function
\begin{equation*}
G(w)=\oint_{\partial D}g(w,z)\,\mathrm{d}z
\end{equation*}
is holomorphic in $\Omega$.
\end{proposition}
\begin{proof}
 In view of Definition \ref{dfanalfunct}, in order to prove holomorhy of $G(w)$, we should prove that it is continuous in $\Omega$ and analytic on affine lines. Denote by $g(w,z)$ restriction of the integrand in $G(w)$ on $\Omega\times\partial D$. Conditions of the proposition imply that it is holomorphic with respect to $w$ in $\Omega$ for each fixed $z\in\partial D$ and continuous in $\Omega\times\partial D$. By Lemma \ref{lmunifcont}, for any $w_0\in\Omega$	and $\epsilon>0$ there is a neighborhood $\mathcal{O}(w_0)$ of $w_0$ such that 
\begin{equation*}
\forall\; w\in\mathcal{O}(w_0):\quad \max_{z\in\partial D} |g(w,z)-g(w_0,z)|<\frac{\epsilon}{length(\partial D)}.
\end{equation*}
Then we obtain easily that $|G(w)-G(w_0)|<\epsilon$. Thus, $G(w)$ is continuous in $\Omega$.  Now consider an affine line $L_{w,\,h}$ of the form \eqref{affline} with $w_0\in\Omega$ and $h\in E$. We have:
\begin{equation*}
G(w_0+\lambda h)=\oint_{\partial D}g(w_0+\lambda h,z)\,\mathrm{d}z.
\end{equation*}	
Thus, we have reached the classical situation of an integral depending on parameter for a function of finite number of scalar complex variables  (\cite{Shab}, Chapt I, \S 2, Lemma in the end of $n^0$ 5). Then the function $G(w_0+\lambda h)$ is holomorphic with respect to $\lambda$ in the domain 
\begin{equation*}
\{\lambda\in\C:\,w_0+\lambda h\in L_{w,\,h}\cap\Omega \}.\qedhere
\end{equation*}
\end{proof}

In the proof of Proposition \ref{pronintpar} we have used the following claim:

\begin{lemma}\label{lmunifcont}
	Let $\Omega$, $E$, $D$ and $\partial D$ be the same as in formulation of Proposition \ref{pronintpar}. Furthermore, let 
	$g(w,z)$ be a function defined in $\Omega\times \partial D$ and continuous there. Then  $g(w,z)$ is continuous at any point $w_0\in\Omega$, uniformly with respect to $z\in\partial D$.
\end{lemma}
\begin{proof}
	Our arguments will be similar to ones in the proof of Theorem 6 from \cite{Tay}. Let $d(\cdot,\,\cdot)$ be a metric, generating topology in $E$.
	Suppose that the claim is false. Then there will exist a number $\epsilon>0$,
	elements $w_k\in\Omega$, and elements $z_k\in\partial D$ such that the inequalities
	\begin{equation*}
	d(w_k,\,w_0)<1/k, \quad |g(w_k, z_k)-g(w_0, z_k)|>\epsilon
	\end{equation*}
	are valid when $k\in\N$. Since $\partial D$ is compact, we may suppose that the points $z_k$ converge to a point $z_0\in\partial D$. Then the inequality
	\begin{equation*}
	\epsilon<|g(w_k, z_k)-g(w_0, z_k)|\le |g(w_k, z_k)-g(w_0, z_0)|+|g(w_0, z_0)-g(w_0, z_k)|
	\end{equation*}
	holds for $k\in\N$. But by the continuity of $g(w,z)$ the right member tends to
	zero with $1/k$, and we are led to a contradiction.
\end{proof}

In the proof of Proposition \ref{prreprWeir} we have used the following claim:

\begin{proposition}\label{prholint}
Let $\Omega$, $E$, $D$ and $\partial D$ be the same as in formulation of Proposition \ref{pronintpar}. 	Let $\Theta(w,\,z)$ and $F(w,\,z)$ be functions, defined and holomorphic in $\Omega\times G$, where $G$ is an open subset of $\C$, containing $\bar D$.
Furthermore, assume that for each fixed $w\in\Omega$ the function $F(w,\,z)$ has $m$ zeros in $D$: $z_1(w),\,z_2(w),\,\dots z_m(w)$ (counting their multiplicities) and it does not vanish on $\partial D$. Then the representation
\begin{equation}\label{genargpin}
\sum_{\nu=1}^m\Theta(w,\,z_\nu(w))=\frac{1}{2\pi i}\oint_{\partial D_R}\Theta(w,\,z)\frac{\partial F(w,z)/\partial z}{F(w,z)}\,\mathrm{d}z.
\end{equation}
is valid, where the integral is a holomorphic function in $\Omega$.
\end{proposition}
\begin{proof}
Notice that for any fixed $w\in\Omega$ representation \eqref{genargpin} follows from Cauchy's Residue Theorem, since as it is not difficult to show, residue of the integrand at its pole $z_\nu(w)$ is equal to $m_\nu\Theta(w,\,z_\nu(w))$, where $m_\nu$ is the multiplicity of the zero $z_\nu(w)$ of $F(w,\,z)$.
By Proposition \ref{analderiv}, the function $\partial F(w,z)/\partial z$ is holomorphic in $\Omega\times G$. Furthermore, by claims (ii) and (iv) of Proposition \ref{properholfunc}, the function $(F(w,z))^{-1}$ is continuous in $\Omega\times\partial D$ and holomorphic in $\Omega$ for any fixed $z\in\partial D$. Then in view of claims (i) and (iii) of Proposition \ref{properholfunc}, the function
\begin{equation*}
g(w,z)=\Theta(w,\,z)\frac{\partial F(w,z)/\partial z}{F(w,z)}
\end{equation*}
satisfies the conditions of Proposition \ref{pronintpar}. Hence the integral in \eqref{genargpin} is a holomorphic function in $\Omega$.
\end{proof}

We used also the following claim:

\begin{proposition}\label{properholfunc}
Let $\Omega$ be an open subset  of a complex Frechet space $E$ and $f:\,\Omega\rightarrow\C$, $g:\,\Omega\rightarrow\C$ be functions. 
\vskip2mm

\noindent (i) If $f(w)$ and $g(w)$ are continuous, then so is their product $f(w)g(w)$;
\vskip2mm

\noindent (ii) If $f(w)$ is continuous and not vanish in $\Omega$, then so is $\frac{1}{f(w)}$;
\vskip2mm

\noindent (iii) If $f(w)$ and $g(w)$ are holomorphic, then so is their product $f(w)g(w)$;
\vskip2mm

\noindent (iv) If $f(w)$ is holomorphic and not vanish in $\Omega$, then so is $\frac{1}{f(w)}$.
 \end{proposition}
\begin{proof}
 Let us take $w_0\in\Omega$ and arbitrary $\epsilon>0$.
 
 (i) We can choose a neighborhood $\mathcal{O}(w_0)$ of $w_0$ so that for any $w\in\mathcal{O}(w_0)$ 
 $|g(w)|<|g(w_0)|+1$, $|f(w)-f(w_0)|<\frac{\epsilon}{2(|g(w_0)|+1)}$ and $|g(w)-g(w_0)|<\frac{\epsilon}{2(|f(w_0)|+1)}$. Then
 \begin{eqnarray*}
 &&|f(w)g(w)-f(w_0)g(w_0|\le\\
 &&|f(w)-f(w_0)||g(w)|+|g(w)-g(w_0)|f(w_0)|<\epsilon.
 \end{eqnarray*}
Claim (i) is proven.

(ii) Let us choose $\mathcal{O}(w_0)$ so that for any $w\in\mathcal{O}(w_0)$ $|f(w)|\ge|f(w_0)|/2|$ and $|f(w)-f(w_0)|<\frac{\epsilon|f(w_0)|^2}{2}$. Then
\begin{equation*}
|\frac{1}{f(w)}-\frac{1}{f(w_0)}|=\frac{|f(w)-f(w_0)|}{|f(w)||f(w_0)|}<\epsilon.
\end{equation*}
Claim (ii) is proven.

(iii) In view of Definition \ref{dfanalfunct} and claim (i), it is enough to show that $f(w)g(w)$ is analytic on any affine line. But this follows from the fact that product of two analytic functions of one complex variable is analytic. Claim (iii) is proven.

Claim (iv) is proved analogously.
\end{proof}

The following claim on hlomorphy of derivative is known:

\begin{proposition}\label{analderiv}[\cite{BS}, Proposition 6.4]
Let $U$ be an open subset in a complex Frechet space $E$. If $f:\,U\rightarrow\C$ is analytic, the function
\begin{equation*}
U\times E\ni(w,h)\rightarrow\delta_w f(h)\in \C
\end{equation*}
is analytic, where $\delta_w f(h)=\frac{\mathrm{d}}{\mathrm{d}\lambda}f(w+\lambda h)_{\lambda=0}$.
\end{proposition}

 We used the following generalization of Hartogs Theorem:

\begin{proposition}\label{prHartogs}[\cite{BS}, Corollary 6.2] 
Let $E_1$ and $E_2$ be complex Frechet spaces and $U$ be an open subset of $E_1\times E_2$.  $f:\;U\ni(x_1,\,x_2)\rightarrow f(x_1,x_2)\in\C$ is separately analytic, then $f$ is analytic.
\end{proposition}

We used also the following known generalization of uniqueness theorem for analytic functions:

\begin{proposition}\label{pruniq}[\cite{BS}, Proposition 6.6, claim I]
Let $U$ be an open connected subset in a complex Frechet space $E$. If $f:\,U\rightarrow\C$ is analytic and $f(w)=0$ in an open subset $V\subset U$, then $f(w)\equiv 0$.
\end{proposition}

{abst}
\end{document}